\documentclass[10.5pt,a4paper]{article}

\usepackage{amssymb}
\usepackage{amsmath}
\usepackage{tensor}
\usepackage{color}
\newtheorem{Th}{Theorem}[section]
\newtheorem{Prop}{Proposition}[section]
\newtheorem{Lm}{Lemma}[section]
\newtheorem{Lma}{Lemma}[section]
\newtheorem{Co}{Corollary}

\newcommand{\be}{\begin{equation}}
\newcommand{\ee}{\end{equation}}
\newcommand{\bes}{\begin{equation*}}
\newcommand{\ees}{\end{equation*}}
\newcommand{\R}{\mathbb{R}}

\newcommand\res{\mathop{\hbox{\vrule height 7pt width .5pt depth 0pt
\vrule height .5pt width 6pt depth 0pt}}\nolimits}

\def\theequation{\thesection.\arabic{equation}}
\def\theTh{\Roman{section}.\arabic{Th}}
\def\theProp{\Roman{section}.\arabic{Prop}}
\def\theCo{\Roman{section}.\arabic{Co}}
\def\theLm{\Roman{section}.\arabic{Lma}}

\def\theRm{\Roman{section}.\arabic{Rm}}
\newcommand{\reset}{\setcounter{equation}{0}\setcounter{Th}{0}\setcounter{Prop}{0}\setcounter{Co}{0}\setcounter{Lma}{0}\setcounter{Rm}{0}}

\def\al{\alpha}
\def\la{\lambda}
\def\eps{\varepsilon}

\def\pro{\pi_{\vec{n}}}
\def\bn{\vec{n}}

\def\bH{\vec{H}}
\def\bL{\vec{L}}
\def\bR{\vec{R}}

\def\bX{\vec{X}}
\def\bU{\vec{U}}
\def\bv{\vec{v}}
\def\bV{\vec{V}}
\def\bA{\vec{A}}

\def\bp{\vec{\Phi}}
\def\bT{\vec{T}}

\def\p{\partial}
\def\di{D_1(0)}

\def\bul{\bullet}

\def\uow{Institute for Mathematics and Its Applications, University of Wollongong, 2522 New South Wales, Australia}

\def\res{\mathop{\hbox{\vrule height 7pt width .5pt 
depth 0pt\vrule height .5pt width 6pt depth 0pt}}\nolimits}

\begin{document}
\title{Analysis of the 
Inhomogeneous Willmore Equation}
\author{Yann Bernard\footnote{School of Mathematical Sciences, Monash University, 3800 Victoria, Australia}, Glen Wheeler\footnote{\uow}, Valentina-Mira Wheeler\footnote{\uow}}
\date{\today}
\maketitle

{\bf Abstract:} {\it We study a class of fourth-order geometric problems
modelling  Willmore surfaces, conformally constrained Willmore surfaces,
isoperimetrically constrained Willmore surfaces, bi-harmonic surfaces in the
sense of Chen, among others. We prove several local energy estimates and
 derive a global gap lemma.}


\reset

\section{Introduction and Main Results}\label{intro}

Let $\Sigma$ be a smooth two-dimensional closed oriented manifold, and let
$g_0$ be a smooth reference metric on $\Sigma$. For any $s\ge1$, the Sobolev
space $W^{k,p}(\Sigma,\R^s)$ is the space of measurable maps
$f:\Sigma\rightarrow\R^s$ for which
\bes
\sum_{j=0}^{k}\,\int_\Sigma|\nabla^jf|^p_{g_0}\,d\text{vol}_{g_0}\;<\;\infty\:.
\ees
For a closed surface $\Sigma$, this space is independent of the reference
metric $g_0$. \\
The notion of {\it weak immersion with $L^2$-bounded second fundamental form}
is well-understood and has been extensively studied (the interested reader will
find a detailed account in \cite{ParksCity} and the references therein). They
will be the main object of study in this paper, and we now recall the main
definition. Let $\bp:\Sigma\rightarrow\R^m$, for $m\ge3$, be measurable and
Lipschitz.
The associated pull-back metric $g:=\bp^*g_{\R^m}$ is given almost everywhere
by
\bes
g(X,Y)\;:=\;d\bp(X)\cdot d\bp(Y)\:,\qquad\forall\:\:X, Y\in T\Sigma\:,
\ees
where dot indicates the standard scalar product in $\R^m$. Unless otherwise
specified, we will demand that $g$ be non-degenerate, that is that there exists
a constant $c>0$ satisfying
\be\label{condi2}
c^{-1}g_0(X,X)\;\leq\;g(X,X)\;\leq\;c\,g_0(X,X)\:,\qquad\forall\:\:X\in T\Sigma\:.
\ee
This makes $(\Sigma, \bp^*g_{\R^m})$ a Riemannian 2-manifold with a rough metric.
The Gauss map is a bounded measurable map $\bn$ taking values in the
Grassmanian $Gr_{m-2}(\R^m)$ of oriented $(m-2)$-planes in $\R^m$ satisfying
\bes
\bn\;:=\;\star\dfrac{\partial_{x^1}\bp\wedge\partial_{x^2}\bp}{|\partial_{x^1}\bp\wedge\partial_{x^2}\bp|}\:,
\ees
where $\star$ denotes the standard Hodge star operator, and $\{x^1,x^2\}$ is an
arbitrary choice of local coordinates. \\
Finally, to say that the weak immersion $\bp$ has square integrable second
fundamental form amounts to requiring that
\be\label{condi3}
\int_\Sigma |d\bn|^2_{g}\,d\text{vol}_g\;<\;\infty\:.
\ee
We let
\bes
\mathcal{E}_{\Sigma}\;:=\;\big\{\bp:\Sigma\rightarrow\R^m\:\:\:\text{measurable and Lipschitz such that (\ref{condi2}) and (\ref{condi3}) hold}\big\}\:.
\ees

Rescaling if necessary, condition (\ref{condi3}) ensures that on some local patch, let us say it is the unit-disk $D_1(0)$, there holds
\be\label{condi4}
\int_{D_1(0)} |\nabla\bn|^2\,dx^1dx^2\;<\;\dfrac{8\pi}{3}\:.
\ee
Here $\{x^1,x^2\}$ are local coordinates on $D_1(0)$ and $\nabla$ stands for the usual flat gradient in these coordinates. A well-known result (\cite{MS, Riv2}) states that if $\bp\in\mathcal{E}_{D_1(0)}$ satisfies (\ref{condi4}), then there exists a bi-Lipschitz homeomorphism $\psi$ of $D_1(0)$ such that the map $\bp\circ\psi:D_1(0)\rightarrow\R^m$ is conformal, namely
\bes
\partial_{x^i}(\bp\circ\psi)\cdot\partial_{x^j}(\bp\circ\psi)\;=\;\text{e}^{2\la}\delta_{ij}\:,
\ees
for some conformal factor $\la$. Without loss of generality, as we are only
concerned with locally analysing the solutions to problems that are independent
of parametrisation, we will henceforth suppose that $\bp$ itself is conformal.\\


The present paper is concerned with studying the local analytical properties of the inhomogeneous Willmore equation. To an immersion $\bp\in\mathcal{E}_{\Sigma}$ of
an oriented two-dimensional manifold $\Sigma$ into $\R^m$, some $m\ge3$, we
assign the second fundamental form $\bA:=\pro D^2\bp$, where $\pro$ denotes the
projection of vectors in $\R^m$ onto the $(m-2)$-place defined by the Gauss map
$\bn$. The trace of the 2-tensor $\bA$ with respect to $g$ is twice the
normal-valued mean curvature vector:
\bes
\bH\;:=\;\dfrac{1}{2}\,\text{Tr}_g\bA\:.
\ees
Willmore immersions are critical points of the Willmore energy
\bes
\int_\Sigma|\bH|^2\,d\text{vol}_g\:.
\ees
The study of Willmore immersions has been steadily gaining momentum over the
last century.
It would be impossible to give a detailed account of the various works and
results that have appeared in recent years.
We content ourselves with mentioning the tour de force by Marques and Neves in
\cite{MN}, where they prove the celebrated Willmore conjecture \cite{Wil}: the
Clifford torus minimizes, up to M\"obius transformations, the Willmore energy
in the class of immersed tori in $\R^3$.
Although the Willmore conjecture is now resolved, the study of Willmore immersions
continues to grow in intensity.
\\

Any critical point of the Willmore energy satisfies the following fourth-order, quasi-linear, strongly coupled system of equations \cite{Wil, Wei}:
\be\label{willeq}
\Delta_\perp\bH+\big\langle\bA\cdot\bH,\bA\big\rangle_g-2|\bH|^2\bH\;=\;\vec{0}\:,
\ee
where $\Delta_\perp$ is the negative covariant Laplacian for the connection in
the normal bundle. The dot indicates the standard scalar product of vectors in
$\R^m$, while the product $\langle\cdot,\cdot\rangle_g$ is the usual
contraction product with respect to the metric $g$ for tensors. Naturally, when
constraints are imposed on the problem of varying the Willmore energy, the
right-hand side of (\ref{willeq}) is no longer zero. Various examples are
provided in \cite{Ber2} and we will below look closer at a few specific cases
of relevance in applications. Thus we are motivated to study a problem of the
type
\be\label{oureq}
\Delta_\perp\bH+\big\langle\bA\cdot\bH,\bA\big\rangle_g-2|\bH|^2\bH\;=\;\vec{\mathcal{W}}\:,
\ee
where the right-hand side $\vec{\mathcal{W}}$ is assumed to be known.
Naturally, $\vec{\mathcal{W}}$ has to be normal vector for (\ref{oureq}) to
make sense. It also has to be independent of parametrization. Before going any
further, an important observation is in order. When $\bp$ lies in
$\mathcal{E}_\Sigma$, it is clear that $\bH$ is square integrable. Even in the
case when $\vec{\mathcal{W}}\equiv\vec{0}$, the term $|\bH|^2\bH$ is already
problematic, for it lies in no space that enables us to understand the equation
in a distributional sense to the equation.
Nevertheless, one may study the problem and obtain estimates, as is done for
example in \cite{Whe1} and the references therein.
Another approach was originally devised by Tristan Rivi\`ere in \cite{Riv1}.
It relies mainly on the fact that the left-hand side of (\ref{oureq}) can
be factored into an exact divergence, thereby rendering possible the
assignment of a distributional sense to \eqref{willeq}.
In \cite{Ber2}, it is shown that the divergence structure seemingly hidden in
(\ref{willeq}) is a direct consequence of Noether's theorem applied to
the translation invariance of the Willmore energy. The present paper should be
understood as a companion to \cite{Ber2}. While in the latter only identities
were derived, the present work brings to fruition the reformulations presented
in \cite{Ber2} by obtaining local analytical results for problems of the type
(\ref{oureq}). The present paper should also be seen as a companion to
\cite{Whe1}, where only a specific class of right-hand sides $\vec{\mathcal{W}}$
were considered. The class of possible right-hand sides will be here
significantly expanded. \\

As was shown in \cite{Riv1}, any conformal immersion
$\bp:D_1(0)\rightarrow\R^m$ that satisfies the Willmore equation (\ref{willeq})
also satisfies the equation
\bes
\text{div}\big(\nabla\bH-2\pro\nabla\bH+|\bH|^2\nabla\bp\big)\;=\;\vec{0}\qquad\text{on}\:\:D_1(0)\:,
\ees
where $\pro$ denotes projection on the normal bundle. The operators $\nabla$
and $\text{div}$ are understood in local coordinates $\{x^1,x^2\}$ on the unit
disk $D_1(0)$.  This motivates us to consider inhomogeneous Willmore problems of
the type
\be\label{ourdiw}
\text{div}\big(\nabla\bH-2\pro\nabla\bH+|\bH|^2\nabla\bp+\bT\big)\;=\;\vec{v}\:,
\ee
for some vector field $\bT\in\Gamma(\R^2\otimes\R^m)$ and some normal vector
field $\bv\in\Gamma(\R^m)$. Many known classes of immersed surfaces satisfy a
problem of this type. 
\begin{itemize}
\item[(1)] Willmore immersions with $\bv\equiv\text{0}$ and $\bT\equiv\vec{0}$.
\item[(2)] Constrained Willmore immersions. 
\begin{itemize}
\item[(i)] Varying the Willmore energy $\int_\Sigma|\bH|^2d\text{vol}_g$ in a
fixed conformal class (i.e. with infinitesimal, smooth, compactly supported,
conformal variations) gives rise to a more general class of surfaces called
{\it conformally-constrained Willmore surfaces} whose corresponding
Euler-Lagrange equation \cite{BPP, KS3, Sch} is expressed as follows. 
Let $\vec{h}_0$ denote the trace-free part of the second fundamental form, namely
\bes
\vec{h}_0\;:=\;\vec{h}\,-\,\bH g\:.
\ees
A conformally-constrained Willmore immersion $\bp$ satisfies
\be\label{cwe}
\Delta_\perp\bH+(\bH\cdot\vec{h}_j^i)\vec{h}^j_i-2|\bH|^2\bH\;=\;\big(\vec{h}_0\big)_{ij}q^{ij}\:,
\ee
where $q$ is a transverse\footnote{i.e. $q$ is divergence-free: $\nabla^j
q_{ji}=0\:\:\forall i$.}  traceless symmetric 2-form. This tensor $q$ plays the
role of Lagrange multiplier in the constrained variational problem. It is shown
in \cite{Ber1, Ber2} that in a conformal parametrization, with conformal
parameter $\la$, (\ref{cwe}) can be brought in the form (\ref{ourdiw}) by
setting $\bv\equiv\vec{0}$ and 
\bes
\bT\;=\;-\,\text{e}^{-2\la}M_q\nabla^\perp\bp\:,
\ees 
where $\nabla^\perp\bp:=(-\partial_{x^2}\bp,\partial_{x^1}\bp)$, and $M_q$ is the matrix
\bes
M_q\;:=\;\left(\begin{array}{cc}-q_{12}&q_{11}\\q_{11}&q_{12}   \end{array}\right)\:.
\ees
\item[(ii)] Bilayer models \cite{BWW, Can, Hef}. These models bear also the name of Helfrich, Canham-Helfrich, and arise in the modelling of the surface of liposomes and vesicles (see \cite{Ber2} and the references therein). One seeks to minimize the Willmore energy under the requirement that the area $A(\Sigma)$, the volume $V(\Sigma)$, and the total curvature
\bes
M(\Sigma)\;:=\;\int_{\Sigma} H\,d\text{vol}_g
\ees
be prescribed. This leads to an equation of the type (\ref{oureq}) with
\bes
\vec{\mathcal{W}}\;=\;2\big(\beta+\alpha H+\gamma K\big)\bn\:,
\ees
where $K$ is the Gauss curvature, and $\al$, $\beta$, $\gamma$ are three given parameters acting as Langrange multipliers.\\
As shown in \cite{Ber2}, this problem can be brought in the form (\ref{ourdiw})
with $\bv\equiv\vec{0}$ and $|\bT|\lesssim1+|\nabla\bn|$. 
\item[(iii)] Another instance in which minimizing the Willmore energy arises is
the isoperimetric problem \cite{KMR, Scy}, which consists in minimizing the
Willmore energy under the constraint that the dimensionless isoperimetric ratio
$\sigma:=36\pi V^2/A^3$ be a given constant in $(0,1]$. As both the Willmore energy
and the constraint are invariant under dilation, one might fix the volume
$V=1/(6\sqrt\pi)$, forcing the area to satisfy $A=\sigma^{1/3}$. This problem is thus
equivalent to the bilayer model with $\gamma=0$ (no constraint imposed on the
total curvature, but the volume and area are prescribed separately).  
\end{itemize}
\item[(3)] Chen surfaces. An isometric immersion $\bp:N^{n}\rightarrow\R^{m>n}$
of an $n$-dimensional Riemannian manifold $N^n$ into Euclidean space is called
{\it biharmonic} if the corresponding mean-curvature vector $\bH$ satisfies 
\be\label{chen}
\Delta_g\bH\;=\;\vec{0}\:.
\ee
The study of biharmonic submanifolds was initiated by B.-Y. Chen \cite{BYC1} in
the mid 1980s as he was seeking a classification of the finite-type
submanifolds in Euclidean spaces. Independently, G.Y. Jiang \cite{Jia} also
studied (\ref{chen}) in the context of the variational analysis of the
biharmonic energy in the sense of Eells and Lemaire. Chen conjectures that a
biharmonic immersion is necessarily minimal\footnote{The conjecture as
originally stated is rather analytically vague: no particular hypotheses on the
regularity of the immersion are {\it a priori} imposed. Many authors consider
only smooth immersions.}.
Smooth solutions of (\ref{chen}) are known to be minimal for $n=1$ \cite{Dim1},
for $(n,m)=(2,3)$ \cite{Dim2}, and for $(n,m)=(3,4)$ \cite{HV}. In \cite{Whe3},
it is shown that Chen's conjecture holds up to a growth condition on the
Willmore energy, and in \cite{BWW2} the parabolic flow with velocity given by Chen's operator is studied.  
Chen's conjecture has been solved under a variety of
hypotheses (see the recent survey paper \cite{BYC2}). The statement remains
nevertheless open in general, and in particular for immersed surfaces in
$\R^m$. In \cite{Ber2b}\footnote{This paper is the precursor to the published
version \cite{Ber2}, which, to the referee's request, no longer
addresses the question of Chen immersions.}, it is shown that Chen surfaces
satisfy an equation of the type (\ref{oureq}) with 
\bes
|\vec{\mathcal{W}}|\;\simeq\;|\bA|^3\;.
\ees
It can more precisely be brought in the form (\ref{ourdiw}) with
$\vec{v}\equiv\vec{0}$ and $|\bT|\lesssim\text{e}^{\la}|\nabla\bn|^2$. 
\item[(4)] Complete Willmore immersions in asymptotically flat spaces also
satisfy a problem of the type (\ref{ourdiw}). Details may be found in
\cite{BR3}.
\item[(5)] Equilibria of flow equations. In \cite{KS1}, stability of the
sphere is proven for the Willmore flow. Global existence is obtained by
contradiction: one assumes that existence time is finite, and then rescales
around a point in space-time where the energy concentrates.  Local estimates
allow one to construct a blowup.  The blowup is shown to be an entire Willmore
surface with small energy.  To this blowup one applies a gap lemma, that
implies any such surface is a standard flat plane.  This is in contradiction
with the concentration of energy hypothesis, and so no such concentration
points can occur, and the flow exists for all time.
This argument is by now standard, having been adapted at least to constrained surface
diffusion flows \cite{Whe,Whe2}, locally constrained Willmore flow \cite{MW2},
Willmore flow in Riemannian spaces \cite{Lnk,MWW}, and a geometric triharmonic
heat flow \cite{MPW}.

An appropriate gap lemma combined with local regularity is crucial and so far
has been established separately for each of the flows given above.
As our work here holds for more general equations than what is currently
available, we expect that the results in this paper will apply to a broad class
of fourth-order evolution equations.
It is an interesting open question to investigate higher-order cases.
\end{itemize}

Our first main result consists of local energy estimates.
\begin{Th}\label{estim}
Let $\bp\in W^{2,2}\cap W^{1,\infty}(D_1(0),\R^m)$ be a conformal immersion with conformal parameter $\la$ satisfying
$$
\Vert\nabla\la\Vert_{L^{2,\infty}(D_1(0))}\;<\;+\infty\:,
$$
where $L^{2,\infty}$ denotes the weak-$L^2$ Marcinkiewicz space. Suppose that 
\be\label{petite}
\int_{D_1(0)}|\nabla\bn|^2\,dx\,=\,\eps_0^2\:.
\ee
Provided that $\eps_0$ is sufficiently small, there is a universal constant $C(\eps_0,\Vert\nabla\la\Vert_{L^{2,\infty}(D_1(0))})$ for which the following statements hold. 
\begin{itemize}
\item[(i)] Let $p\in(1,\infty)$. Suppose that $\bp$ is a solution on $D_1(0)$ of
\bes
\text{div}\big(\nabla\bH-2\pro\nabla\bH+|\bH|^2\nabla\bp+\bT\big)\;=\;\vec{0}\:.
\ees
Then for all $D_\rho(x)\subset D_1(0)$, we have
\bes
\rho^{2-\frac{2}{p}}\Vert\nabla^2\bn \Vert_{L^p(D_{\rho/2}(x))}\;\leq\;(M+1)^2\:,
\ees
where
\bes
M\;:=\;C(\eps_0,\Vert\nabla\la\Vert_{L^{2,\infty}(D_1(0))})\Big[\rho^{2-\frac{2}{p}}\Vert\text{e}^\la\bT\Vert_{L^p(D_\rho(x))}+\Vert\nabla\bn\Vert_{L^2(D_\rho(x))}\Big]\:.
\ees
Moreover,
\begin{eqnarray*}
\left\{\begin{array}{lclcl}
\rho^{2-\frac{2}{p}}\Vert\nabla\bn \Vert_{L^{2p/(2-p)}(D_{\rho/2}(x))}&\leq&(M+1)^2&,&\text{if}\:\:p\in(1,2)\\[1.5ex]
\rho^{1-\frac{2}{q}}\Vert\nabla\bn \Vert_{L^{q}(D_{\rho/2}(x))}&\leq&M+1&,&\forall\:\:q<\infty\:,\:\:\text{if}\:\:p=2\\[1.5ex]
\rho\Vert\nabla\bn \Vert_{L^{\infty}(D_{\rho/2}(x))}&\leq&M+1&,&\text{if}\:\:p>2\:,
\end{array}\right.
\end{eqnarray*}
and
\begin{eqnarray*}
\left\{\begin{array}{lclcl}
\rho^{2-\frac{2}{p}}\Vert\text{e}^{\la}\bH \Vert_{L^{2p/(2-p)}(D_{\rho/2}(x))}&\leq&M(M+1)&,&\text{if}\:\:p\in(1,2)\\[1.5ex]
\rho^{1-\frac{2}{q}}\Vert\text{e}^{\la}\bH \Vert_{L^{q}(D_{\rho/2}(x))}&\leq&M&,&\forall\:\:q<\infty\:,\:\:\text{if}\:\:p=2\\[1.5ex]
\rho\Vert\text{e}^{\la}\bH \Vert_{L^{\infty}(D_{\rho/2}(x))}&\leq&M&,&\text{if}\:\:p>2\:.
\end{array}\right.
\end{eqnarray*}

\item[(ii)] Let $r\in[1,\infty)$. Suppose that $\bp$ is a solution on $D_1(0)$ of
\bes
\text{div}\big(\nabla\bH-2\pro\nabla\bH+|\bH|^2\nabla\bp\big)\;=\;\vec{v}\:,
\ees
with $\text{e}^{\la}\bv\in L^r(\di)$. For all $D_\rho(x)\subset D_1(0)$, we have
\bes
\left\{\begin{array}{rclcl}
\rho^{2-\frac{2}{p}}\Vert\nabla^2\bn\Vert_{L^p(D_{\rho/2}(x))}+\rho^{2-\frac{2}{p}}\Vert\nabla\bn\Vert_{L^{2p/(2-p)}(D_{\rho/2}(x))}&\leq&(M+1)^2&,&\quad\forall\:\:p\in(1,2)\:\:\:\text{if}\:\:\:r=1\\[2.5ex]
\rho^{3-\frac{2}{r}}\Vert\nabla^2\bn\Vert_{L^{2r/(2-r)}(D_{\rho/2}(x))}&\leq&(M+1)^2&,&\quad\text{if}\:\:\:r\in(1,2)\\[2.5ex]
\rho^{2-\frac{2}{q}}\Vert\nabla^2\bn\Vert_{L^{q}(D_{\rho/2}(x))}&\leq&(M+1)^2&,&\quad\forall\:\:q<\infty\:\:\:\text{if}\:\:\:r=2\:,\\[2.5ex]
\rho^2\Vert\nabla^2\bn\Vert_{L^{\infty}(D_{\rho/2}(x))}&\leq&(M+1)^2&,&\quad\text{if}\:\:\:r>2\:,
\end{array}\right.
\ees
with
\bes
M\;=\;C(\eps_0,\Vert\nabla\la\Vert_{L^{2,\infty}(D_1(0))})\Big[\rho^{3-\frac{2}{r}}\Vert\text{e}^{\la}\bv\Vert_{L^{r}(D_\rho(x))}+\Vert\nabla\bn\Vert_{L^2(D_\rho(x))}\Big]\:.
\ees
Furthermore, if $r>1$, the following estimates hold:
\bes
\rho\Vert\nabla\bn\Vert_{L^\infty(D_{\rho/2}(x))}\;\leq\;M+1\:.
\ees
and
\bes
\rho^{3-\frac{2}{r}}\Vert\nabla^3\bn\Vert_{L^r(D_{\rho/3}(x))}\;\leq\;(M+1)^3\:.
\ees
\end{itemize}
\end{Th}

Theorem \ref{estim} is used to prove the following regularity result:
\begin{Co}\label{smooth}
Let $\bp\in W^{2,2}\cap W^{1,\infty}(D_1(0),\R^m)$ be a conformal immersion satisfying (\ref{ourdiw}) on the disk $D_1(0)$. If $\bT$ and $\bv$ are smooth, so is $\bp$. 
\end{Co}

Finally, we derive an interesting geometric ``gap" result, obtained using the same techniques as those leading to Theorem \ref{estim}. 

\begin{Th}\label{gap0}
Let $\Sigma$ be a connected oriented complete immersed surface in $\R^m$ whose mean curvature vector satisfies an inhomogeneous Willmore problem of the type\footnote{We use the same notation as in (\ref{willeq}).}
\bes
\Delta_\perp\bH+\big\langle\bA\cdot\bH,\bA\big\rangle_g-2|\bH|^2\bH\;=\;\text{O}\big(|\bA|^3\big)\:.
\ees 
There exists an $\eps_0>0$ such that if 
\bes
\int_\Sigma|\bA|^2d\text{vol}_g\,<\,\eps_0^2\:,
\ees
then $\Sigma$ is a flat plane. 
\end{Th}

This gap result is to be compared to the one given in \cite{Whe1} (see also \cite{MW1}). \\

A word of caution is now in order. Should $\bp$ be a (conformal) Willmore
immersion satisfying the small energy condition (\ref{petite}), then
$\bv\equiv\vec{0}$ and Theorem \ref{estim}-(i) gives the estimate 
\bes
\Vert\nabla\bn\Vert_{L^{\infty}(D_{\rho/2}(x))}\;\leq\;C\rho^{-1}\big(1+\Vert\nabla\bn\Vert_{L^2(D_\rho(x))}\big)\:.
\ees
This estimate, which we will term {\it parametric $\eps$-regularity}, is the
one that was originally derived by Rivi\`ere in \cite{Riv1}.
In conformal parametrization, $|\nabla\bn|\simeq\text{e}^{\la}\bA$, where $\bA$ is the second fundamental form, so the above reads
\be\label{compar0}
\Vert\text{e}^{\la}\bA\Vert_{L^{\infty}(D_{\rho/2}(x))}\;\leq\;C\rho^{-1}\big(1+\Vert\text{e}^{\la}\bA\Vert_{L^2(D_\rho(x))}\big)\:.
\ee
Knowing that our conformal immersion does not ``distort'' flat disks much, we can further rephrase \eqref{compar0} as
\be\label{compar1}
\Vert\bA\Vert_{L^{\infty}(D^g_{\rho/2}(x))}\;\leq\;C\rho^{-1}\big(1+\Vert\bA\Vert_{L_g^2(D^g_\rho(x))}\big)\:.
\ee
where $D^g_\rho(x)$ is the metric disk with respect to the induced metric pull-back metric $g=\bp^*g_{\R^m}$, and $L^2_g$ is the space $(L^2,d\text{vol}_g)$. The estimate (\ref{compar1}) is to be compared with Kuwert and Sch\"atzle's original estimate \cite{KS1}, which we will call {\it ambient $\eps$-regularity}, and which states that if $\bp:\Sigma\rightarrow\R^m$ is a Willmore immersion with
\bes
\int_{\bp^{-1}(B_\sigma(p))}|\vec{A}|^2d\text{vol}_g\;<\;\eps_0^2
\ees
for some Euclidean ball $B_\sigma(p)\subset\R^m$, and $\eps_0$ is sufficiently small, then
\be\label{compar2}
\Vert\bA\Vert_{L^\infty(\bp^{-1}(B_{\sigma/2}(p)))}\;\leq C\sigma^{-1}\Vert\bA\Vert_{L^2(\bp^{-1}(B_{\sigma}(p)))}\:.
\ee
Remark 2.11 in \cite{KS1} and more explicitly equation (2.18) in \cite{KS2} state that this estimate implies 
\be\label{compar3}
\Vert\bA\Vert_{L^{\infty}(D^g_{\rho/2}(x))}\;\leq\;C\rho^{-1}\Vert\bA\Vert_{L_g^2(D^g_\rho(x))}\:,
\ee
which is manifestly different from (\ref{compar1}). To the authors' knowledge,
it is unclear that (\ref{compar3}) follows from (\ref{compar2}). The two
versions of the $\eps$-regularity, parametric and ambient, are resolutely
distinct and we do not know how to recover one from the other. We suspect that
(\ref{compar3}) might in fact be false, although what does remain true, as
given in Theorem \ref{estim}-(i), is
\bes
\Vert\text{e}^{\la}\bH\Vert_{L^{\infty}(D_{\rho/2}(x))}\;\leq\;C\rho^{-1}\Vert\nabla\bn\Vert_{L^2(D_\rho(x))}\:.
\ees
Fortunately, this estimate is in fact all which is required to correct the
proof of \cite{KS2} and the end-results remain intact. A similar inconsequent
error is found in \cite{BR1}.


\setcounter{equation}{0} 
\reset

\section{Proofs of the Results}
\reset
\subsection{Controlling the conformal factor}
Using F. H\'elein's method of moving Coulomb frames \cite{Hel}, a weak immersion $\bp\in W^{2,2}_{imm}(D_1(0),{\R}^m)$ of the unit disk $D_1(0)$ into $\R^m$ can be reparametrized by a diffeomorphism of $D_1(0)$ to become conformal. Our problem being independent of parametrization, we will without loss of generality suppose that $\bp$ is conformal with parameter $\la$, namely:
\bes
\partial_{x_i}\bp\cdot\partial_{x_j}\bp\;=\;\text{e}^{2\la}\delta_{ij}\:.
\ees
We will henceforth use the notation $\nabla$, $\text{div}$, and $\Delta$ to denote the usual gradient, divergence, and Laplacian operators in flat local coordinates $\{x_1,x_2\}$. \\
Assume
\[
\int_{D_1(0)}|\nabla\vec{n}|^2\ dx=:\eps_0^2\le 8\pi/3\quad \mbox{ and }\quad \|\nabla\la\|_{L^{2,\infty}(D_1(0))}<+\infty\:.
\]
We can call upon Lemma 5.1.4 in \cite{Hel} to deduce the existence of an orthogonal frame $\{\vec{e}_1,\vec{e}_2\}\in W^{1,2}(D_1(0))$ satisfying $\star\bn=\vec{e}_1\wedge\vec{e}_2$ and
\bes
\Vert\nabla\vec{e}_1\Vert_{L^2(D_2(0))}+\Vert\nabla\vec{e}_2\Vert_{L^2(D_2(0))}\;\leq\;C\, \|\nabla\vec{n}\|_{L^2(D_2(0))}\:,
\ees
 As is easily verified, the conformal parameter satisfies
\bes
\Delta\la\;=\;\nabla\vec{e}_1\cdot\nabla^\perp\vec{e}_2\qquad\text{in}\:\:D_1(0)\:.
\ees
Let $\mu$ satisfy
\bes
\left\{\begin{array}{rclcl}
\Delta\mu&=&\nabla\vec{e}_1\cdot\nabla^\perp\vec{e}_2&,&\text{in}\:\:D_1(0)\\[1ex]
\mu&=&0&,&\text{on}\:\:\partial D_1(0)\:.
\end{array}\right.
\ees
Standard Wente estimates (cf. Theorem 3.4.1 in \cite{Hel}) give
\be\label{wenteagain}
\Vert \mu\Vert_{L^\infty(D_1(0))}+\Vert \nabla\mu\Vert_{L^2(D_1(0))}\;\leq\;\Vert\nabla\vec{e}_1\Vert_{L^2(D_1(0))}\Vert\nabla\vec{e}_2\Vert_{L^2(D_1(0))}\;\leq\;  C\, \|\nabla\vec{n}\|^2_{L^2(D_1(0))}  \:.
\ee
The harmonic function $\nu:=\la-\mu$ satisfies the usual estimate
\bes
\int_{D}|\nu-\bar{\nu}|\,dx\;\leq\;C\,\Vert\nabla\nu\Vert_{L^1(D_1(0))}\;\leq\;C\,\Vert\nabla\nu\Vert_{L^{2,\infty}(D_1(0))}\:,
\ees
where $\bar{\nu}$ denotes the average of $\nu$ on the proper subdisk $D\subset\subset D_1(0)$. Hence
\bes
\Vert\nu-\bar{\nu}\Vert_{L^\infty(D)}\;\leq\;C\,\Vert\nabla\nu\Vert_{L^{2,\infty}(D_1(0))}\:,
\ees 
Combining the latter to (\ref{wenteagain}) yields now
\bes
\Vert\la-\bar{\la}\Vert_{L^\infty(D)}\;\leq\;C\,\Vert\nabla\la\Vert_{L^{2,\infty}(D_1(0))}+\;  C\, \|\nabla\vec{n}\|^2_{L^2(D_1(0))} \;\leq\;C(\eps_0,\|\nabla\la\|_{L^{2,\infty}(D_1(0))})\:,
\ees
where $\bar{\la}$ denotes the average of $\la$ on $D$. We can summarize this subsection
by stating the following lemma.
\begin{Lma}
\label{harnack}
Let $\bp\in W^{2,2}_{imm}(D_1(0),{\R}^m)$ be a conformal weak immersion such that
\[
\int_{D_1(0)}|\nabla\vec{n}|^2\ dx=:\eps_0^2\le 8\pi/3\quad \mbox{ and }\quad \|\nabla\la\|_{L^{2,\infty}(D_1(0))}<+\infty\:,
\]
with $\text{e}^\la:=|\p_{x_1}\vec{\Phi}|=|\p_{x_2}\vec{\Phi}|$. Then the following estimate holds for any proper subdisk $D\subset\subset D_1(0)$:
\be
\label{harnacko}
\|\text{e}^{\la}\|_{L^\infty(D)}\|\text{e}^{-\la}\|_{L^\infty(D)}\;\leq\;C(\eps_0,\|\nabla\la\|_{L^{2,\infty}(D_1(0))}) \:.
\ee
\hfill $\Box$
\end{Lma}

\subsection{Proof of Theorem \ref{estim}-(i)}

Per the discussion in the introduction and our aim to study only local properties of solutions to (\ref{oureq}), we assume without loss of generality that the immersion $\bp$ is conformal, i.e. in local coordinates $\{x^1,x^2\}$ on the unit disk $D_1(0)$ that
\bes
\partial_{x^i}\bp\cdot\partial_{x^j}\bp\;=\;\text{e}^{2\la}\delta_{ij}\:,
\ees
with bounded conformal parameter $\la$, and such that $\text{e}^{\la}$ satisfies the Harnack inequality (\ref{harnacko}). We will first begin by studying an inhomogeneous  Willmore equation of the form
\be\label{ourdiv2}
\text{div}\big(\nabla\bH-2\pro\nabla\bH+|\bH|^2\nabla\bp+\bT  \big)\;=\;\vec{0}\:,\qquad\text{on}\:\:D_1(0)\:,
\ee
where $\bT$ satisfies the following condition for some $p\in(1,\infty)$:
\be\label{hypovv}
\Vert\text{e}^{\la}\bT\Vert_{L^p(D_1(0))}\,<\,\infty\:.
\ee
Let $D_\rho(x)\subset D_1(0)$. As is done in \cite{Ber2}, we consider the following two problems
\be\label{eqX}
\Delta\bX\;=\;\nabla\bp\wedge\bT\qquad\text{and}\qquad\Delta Y\;=\;\nabla\bp\cdot\bT\qquad\text{on}\:\:D_\rho(x)\:
\ee
with boundary conditions $\bX|_{\partial D_\rho(x)}=\vec{0}$ and $Y|_{\partial D_\rho(x)}=0$. Standard Calderon-Zygmund estimates give:
\be\label{meluch1}
\Vert\nabla^2\bX\Vert_{L^{p}(D_\rho(x))}+\Vert\nabla^2Y\Vert_{L^{p}(D_\rho(x))}\;\lesssim\;\Vert\text{e}^{\la}\bT\Vert_{L^p(D_\rho(x))}\:,
\ee
up to a universal multiplicative constant. Hence
\be\label{meluch2}
\Vert\nabla\bX\Vert_{L^{2,\infty}(D_\rho(x))}+\Vert\nabla Y\Vert_{L^{2,\infty}(D_\rho(x))}\;\lesssim\;\rho^{2-\frac{2}{p}}\Vert\text{e}^{\la}\bT\Vert_{L^p(D_\rho(x))}\:.
\ee

We now follow the procedure outlined in \cite{Ber2}. Integrating (\ref{ourdiv2}), we infer the existence of a potential $\bL$ satisfying
\be\label{defL}
\nabla^\perp\bL\;=\;\nabla\bH-2\pro\nabla\bH+|\bH|^2\nabla\bp +\bT \;\equiv\;-\nabla\bH+2\pi_T\nabla\bH+|\bH|^2\nabla\bp+\bT\:,
\ee
where $\pi_T$ is the tangential projection. An elementary computation (c.f. equation (II.6) in \cite{Ber1}) reveals that
\be\label{store}
\big|\pi_T\nabla\bH\big|\;\lesssim\;\text{e}^{\la}|\nabla\bn|^2\:.
\ee
As $\bL$ is defined up to an arbitrary constant, we are certainly free to require that
\bes
\int_{D_\rho(x)}\bL\;=\;\vec{0}\:.
\ees
Observe next that
\bes
\Vert\nabla\bH\Vert_{W^{-1,2}(D_\rho(x))}\;\leq\;\Vert\bH\Vert_{L^2(D_\rho(x))}\;\leq\;\Vert\text{e}^{-\la}\Vert_{L^\infty(D_\rho(x))}\Vert\nabla\bn\Vert_{L^2(D_\rho(x))}\:,
\ees
and, owing to (\ref{hypovv}) and (\ref{store}), 
\bes
\big\Vert\bT+2\pi_T\nabla\bH+|\bH|^2\nabla\bp\big\Vert_{L^1(D_\rho(x))}\;\lesssim\;\Vert\text{e}^{-\la}\Vert_{L^{\infty}(D_\rho(x))}\Big[\Vert\text{e}^{\la}\bT\Vert_{L^1(D_\rho(x))}+\Vert\nabla\bn\Vert^2_{L^2(D_\rho(x))}\Big]
\ees
up to a multiplicative constant independent of the parametrization and of the
mean curvature, and irrelevant to our purpose. Geared with these last
inequalities, we call upon Lemma \ref{L2weak} from the Appendix and conclude
that
\be\label{daban}
\Vert\bL\Vert_{L^{2,\infty}(D_\rho(x))}\;\lesssim\;\Vert\text{e}^{-\la}\Vert_{L^{\infty}(D_\rho(x))}\Big[\Vert\text{e}^{\la}\bT\Vert_{L^1(D_\rho(x))}+\Vert\nabla\bn\Vert^2_{L^2(D_\rho(x))}\Big]\:,
\ee
where $L^{2,\infty}$ is the weak-$L^2$ Marcinkiewicz space, seen here as a Lorentz space \cite{Tar}.
Per Lemma \ref{harnack}, $\text{e}^{\la}$ satisfies a Harnack inequality. The above then yields
\be\label{gtn}
\Vert\text{e}^{\la}\bL\Vert_{L^{2,\infty}(D_\rho(x))}\;\lesssim\;\Vert\text{e}^{\la}\bT\Vert_{L^1(D_\rho(x))}+\Vert\nabla\bn\Vert_{L^2(D_\rho(x))}\:.
\ee
We will use the symbol $\lesssim$ to indicate the presence of a multiplicative constant depending at most only on $\eps_0$ and on $\|\nabla\la\|_{L^{2,\infty}(D_1(0))}$. \\

It is shown in \cite{Ber2} that two important identities hold, namely
\bes
\text{div}\big(\bL\wedge\nabla^\perp\bp+\bH\wedge\nabla\bp+\nabla\bX  \big)\;=\;\vec{0}\qquad\text{and}\qquad \text{div}\big(\bL\cdot\nabla^\perp\bp+\nabla Y  \big)\;=\;0\:.
\ees
Again, we infer the existence of two potentials $\bR$ and $S$ satisfying
\be\label{defS}
\nabla\bR\;=\;\bL\wedge\nabla\bp-\bH\wedge\nabla^\perp\bp-\nabla^\perp\bX   \qquad\text{and}\qquad\nabla S\;=\;\bL\cdot\nabla\bp-\nabla^\perp Y\:.
\ee
Owing to (\ref{gtn}) and (\ref{meluch2}), we find that $\nabla\bR$ and $\nabla S$ lie in the weak space $L^{2,\infty}$, namely
\begin{eqnarray*}
&&\hspace{-1cm}\Vert\nabla\bR\Vert_{L^{2,\infty}(D_\rho(x))}+\Vert\nabla S\Vert_{L^{2,\infty}(D_\rho(x))}\nonumber\\[1ex]
&&\hspace{1cm}\lesssim\:
\Vert\text{e}^{\la}\bT\Vert_{L^1(D_\rho(x))}+\Vert\nabla\bn\Vert^2_{L^2(D_\rho(x))}+\Vert\nabla\bn\Vert_{L^2(D_\rho(x))}
+\Vert\nabla\bX\Vert_{L^{2,\infty}(D_\rho(x))}+\Vert\nabla Y\Vert_{L^{2,\infty}(D_\rho(x))}\nonumber\\[1ex]
&&\hspace{1cm}\lesssim\:\rho^{2-\frac{2}{p}}\Vert\text{e}^{\la}\bT\Vert_{L^p(D_\rho(x))}+\Vert\nabla\bn\Vert_{L^2(D_\rho(x))}\:,
\end{eqnarray*}
where $C(\eps_0)$ is a constant depending only on $\eps_0$. In other words
\be\label{estimRS}
\Vert\nabla\bR\Vert_{L^{2,\infty}(D_\rho(x))}+\Vert\nabla S\Vert_{L^{2,\infty}(D_\rho(x))}\;\lesssim\;M\:,
\ee
where for notational convenience, we have set
\be\label{defM}
M\;:=\;\rho^{2-\frac{2}{p}}\Vert\text{e}^{\la}\bT\Vert_{L^p(D_\rho(x))}+\Vert\nabla\bn\Vert_{L^2(D_\rho(x))}\:.
\ee

It is remarkable that $\bR$ and $S$ are linked together via an interesting system of equations that displays a very particular structural type. It is shown in \cite{Ber2} that\footnote{Refer to the Appendix for the notation.}
\be\label{sysRS}
\left\{\begin{array}{lcl}
\Delta\bR&=&\nabla(\star\bn)\bullet\nabla^\perp\bR+\nabla(\star\bn)\cdot\nabla^\perp S+\text{div}\big((\star\bn)\bullet\nabla\bX+(\star\bn)\nabla Y   \big)\\[1ex]
\Delta S&=&\nabla(\star\bn)\cdot\nabla^\perp\bR+\text{div}\big((\star\bn)\cdot\nabla\bX  \big)\:.
\end{array}\right.
\ee
The apparent notational complication is an artifice of codimension only.\\
The system (\ref{sysRS}) is in divergence form. Owing to $\Vert\bn\Vert_{L^\infty(D_1(0))}=1$, to $\Vert\nabla\bn\Vert_{L^2(D_1(0))}<\eps_0$, and to (\ref{meluch1}), we can call upon Proposition \ref{morreydecay}\footnote{Proposition \ref{morreydecay} is proved for one equation, but it is easily adapted for systems. Details are left to the reader.}, which states that for\footnote{i.e. $s<2/(2-p)$ if $p\in(1,2)$, and $s<\infty$ for $p\ge2$.} all $s<2/(2-p)_+$ it holds
\begin{eqnarray}\label{pjo1}
&&\hspace{-2.5cm}\Vert\nabla\bR \big\Vert_{L^{s}(D_{5\rho/8}(x))}+\Vert\nabla S\Vert_{L^{s}(D_{5\rho/8}(x))}\nonumber\\[1ex]
&&\hspace{1cm}\lesssim\:\rho^{\frac{2}{s}-1}\Big[\Vert\nabla\bR\Vert_{L^{2,\infty}(D_\rho(x))}+\Vert\nabla S\Vert_{L^{2,\infty}(D_\rho(x))}+\rho^{2-\frac{2}{p}}\Vert\text{e}^{\la}\bT\Vert_{L^{p}(D_\rho(x))}\Big]\nonumber\\[1ex]
&&\hspace{1cm}\lesssim\:\rho^{\frac{2}{s}-1}M\:,
\end{eqnarray}
where we have used (\ref{estimRS}). Note that (\ref{pjo1}) holds in particular for $s=2p$.

A useful identity is derived in \cite{Ber2}; it relays information on $\bR$ and $S$ back to the immersion $\bp$, namely:
\be\label{backimm}
\text{e}^{2\la}\bH\;=\;\big(\nabla\bR+\nabla^\perp\bX\big)\bullet\nabla^\perp\bp+\big(\nabla S+\nabla^\perp Y\big)\cdot\nabla^\perp\bp\:.
\ee
It follows from this identity and from (\ref{meluch1}) and (\ref{pjo1}) that
\be\label{simone}
\Vert\text{e}^{2\la}\bH\Vert_{L^{2p}(D_{5\rho/8}(x) )}\;\lesssim\;\rho^{\frac{1}{p}-1}\Vert\text{e}^{\la}\Vert_{L^\infty(D_\rho(x))}M\:,
\ee
where as always the symbol $\:\lesssim\:$ indicates the presence of a multiplicative constant involving at most $\eps_0$ and $\Vert\nabla\la\Vert_{L^{2,\infty}(D_1)}$. 
We now use the equation $\Delta\bp=2\text{e}^{2\la}\bH$ by writing $\bp=\bp_1+\bp_0$, where 
\bes
\left\{\begin{array}{rcl}
\Delta\bp_1&=&2\text{e}^{2\la}\bH\\[1ex]
\bp_1&=&\vec{0}
\end{array}\right.
\qquad\text{and}\qquad
\left\{\begin{array}{rclcl}
\Delta\bp_0&=&\vec{0}&,&\text{in}\:\:D_{5\rho/8}(x)\\[1ex]
\bp_0&=&\bp&,&\text{on}\:\:\partial D_{5\rho/8}(x)\:.
\end{array}\right.
\ees
Estimates for $\bp_1$ follow directly from (\ref{simone}), while $\bp_0$ is handled with the help of standard estimates for harmonic functions. We obtain
\begin{eqnarray}\label{tita0}
\Vert\nabla^2\bp\Vert_{L^{2p}(D_{9\rho/16}(x))}&\lesssim&\rho^{\frac{1}{p}-1}\Vert\nabla\bp\Vert_{L^\infty(D_{5\rho/8}(x))}+\Vert\text{e}^{\la}\Vert_{L^\infty(D_\rho(x))}\rho^{\frac{1}{p}-1}M\nonumber\\[1ex]
&\lesssim&\Vert\text{e}^{\la}\Vert_{L^\infty(D_\rho(x))}\rho^{\frac{1}{p}-1}(1+M)\:.
\end{eqnarray}
Using again that $|\nabla\bn|\leq\text{e}^{-\la}|\nabla^2\bp|$ and (\ref{harnacko}), the latter gives
\be\label{nas0}
\Vert\nabla\bn\Vert_{L^{2p}(D_{9\rho/16}(x))}\;\lesssim\;\rho^{\frac{1}{p}-1}(1+M)\:.
\ee
Next, using (\ref{sysRS}), we have
\begin{eqnarray*}
|\Delta\bR|+|\Delta S|&\lesssim&|\Delta\bX|+|\Delta Y|+|\nabla\bn|\big(|\nabla\bR|+|\nabla S|+|\nabla\bX|+|\nabla Y|  \big)\nonumber\\[1ex]
&\lesssim&|\text{e}^{\la}\bT|+|\nabla\bn|\big(|\nabla\bR|+|\nabla S|+|\nabla\bX|+|\nabla Y|  \big)\:.
\end{eqnarray*}
Hence, from (\ref{nas0}), (\ref{meluch1}), and (\ref{pjo1}),
\begin{eqnarray}\label{nadege0}
&&\hspace{-.8cm}\big\Vert\nabla^2\bR\big\Vert_{L^p(D_{17\rho/32}(x))}+\big\Vert\nabla^2S\big\Vert_{L^p(D_{17\rho/32}(x))}\;\;\lesssim\;\;\rho^{\frac{2}{p}-2}\big(\Vert\nabla R\Vert_{L^{2,\infty}(D_\rho(x))}+\Vert\nabla S\Vert_{L^{2,\infty}(D_\rho(x))} \big)+\Vert\text{e}^{\la}\bT\Vert_{L^p(D_{9\rho/16}(x))}\nonumber\\[1ex]
&&\hspace{-.8cm}\:+\Vert\nabla\bn\Vert_{L^{2p}(D_{9\rho/16}(x))}\big(\Vert\nabla\bR\Vert_{L^{2p}(D_{9\rho/16}(x))}+\Vert\nabla S\Vert_{L^{2p}(D_{9\rho/16}(x))}+\Vert\nabla\bX\Vert_{L^{2p}(D_{9\rho/16}(x))}+\Vert\nabla Y\Vert_{L^{2p}(D_{9\rho/16}(x))}  \big)\nonumber\\[1ex]
&&\hspace{-.8cm}\lesssim\:\Vert\text{e}^{\la}\bT\Vert_{L^p(D_{\rho}(x))}+\rho^{\frac{2}{p}-2}M(M+1)\nonumber\\[1ex]
&&\hspace{-.8cm}\lesssim\:\rho^{\frac{2}{p}-2}M(M+1)\:,
\end{eqnarray}
where we have used (\ref{estimRS}).\\
We can now combine (\ref{meluch1}), (\ref{tita0}), and (\ref{nadege0}) into (\ref{backimm}) to find
\begin{eqnarray}\label{wo1}
&&\Vert\nabla(\text{e}^{2\la}\bH)\Vert_{L^p(D_{17\rho/32}(x))}\nonumber\\[1ex]
&&\lesssim\:\:\Vert\text{e}^{\la}\Vert_{L^\infty(D_\rho(x))}\big(\Vert\nabla^2\bR\Vert_{L^p(D_{17\rho/32}(x))}+\Vert\nabla^2S\Vert_{L^p(D_{17\rho/32}(x))} +\Vert\nabla^2\bX\Vert_{L^p(D_{\rho}(x))}+\Vert\nabla^2Y\Vert_{L^p(D_{\rho}(x))}  \big)\nonumber\\
&&+\:\:\Vert\nabla^2\bp\Vert_{L^{2p}(D_{9\rho/16}(x))}\big(\Vert\nabla\bR\Vert_{L^{2p}(D_{3\rho/4}(x))}+\Vert\nabla S\Vert_{L^{2p}(D_{3\rho/4}(x))} +\Vert\nabla\bX\Vert_{L^{2p}(D_{\rho}(x))}+\Vert\nabla Y\Vert_{L^{2p}(D_{\rho}(x))}  \big)\nonumber\\[1ex]
&&\lesssim\:\:\Vert\text{e}^{\la}\Vert_{L^\infty(D_\rho(x))}\Big[\rho^{\frac{2}{p}-2}M(M+1)+\rho^{\frac{2}{p}-2}(M+1)\big(\Vert\text{e}^{\la}\bT\Vert_{L^p(D_{\rho}(x))}+M  \big)  \Big]\nonumber\\[1ex]
&&\lesssim\:\:\Vert\text{e}^{\la}\Vert_{L^\infty(D_\rho(x))}\rho^{\frac{2}{p}-2}M(M+1)\:,
\end{eqnarray}
where we have also used (\ref{meluch1}) and (\ref{pjo1}). From the latter and (\ref{harnacko}), it follows that
\begin{eqnarray}\label{estimH}
\left\{\begin{array}{lclcl}
\Vert\text{e}^{\la}\bH \Vert_{L^{2p/(2-p)}(D_{\rho/2}(x))}&\lesssim&\rho^{\frac{2}{p}-2}M(M+1)&,&\text{if}\:\:p\in(1,2)\\[1.5ex]
\Vert\text{e}^{\la}\bH \Vert_{L^{q}(D_{\rho/2}(x))}&\lesssim&\rho^{\frac{2}{q}-1}M(M+1)&,&\forall\:\:q<\infty\:,\:\:\text{if}\:\:p=2\\[1.5ex]
\Vert\text{e}^{\la}\bH \Vert_{L^{\infty}(D_{\rho/2}(x))}&\lesssim& \rho^{-1}M(M+1)&,&\text{if}\:\:p>2\:.
\end{array}\right.
\end{eqnarray}

Owing to (\ref{backimm}), we verify easily that
\be\label{rozen2}
2\bH\wedge\nabla^\perp\bp\;=\;\big(\nabla^\perp\bR-\nabla\bX)\bul(\star\bn)+\big(\nabla^\perp S-\nabla Y)(\star\bn)\:,
\ee
so that, using (\ref{eqX}) and (\ref{sysRS}),
\begin{eqnarray*}
\big|\text{div}\big(\bH\wedge\nabla^\perp\bp\big)\big|&\lesssim&|\Delta\bX|+|\Delta Y|+|\nabla\bn|\big(|\nabla\bR|+|\nabla S|+|\nabla\bX|+|\nabla Y|  \big)\nonumber\\[1ex]
&\lesssim&|\text{e}^{\la}\bT|+|\nabla\bn|\big(|\nabla\bR|+|\nabla S|+|\nabla\bX|+|\nabla Y|  \big)\:.
\end{eqnarray*}
Hence, from (\ref{nas0}), (\ref{meluch1}), and (\ref{pjo1}),
\begin{eqnarray}\label{nadege}
&&\hspace{-.8cm}\big\Vert\text{div}\big(\bH\wedge\nabla^\perp\bp\big)\big\Vert_{L^p(D_{9\rho/16}(x))}\;\;\lesssim\;\;\Vert\text{e}^{\la}\bT\Vert_{L^p(D_{9\rho/16}(x))}\nonumber\\[1ex]
&&\hspace{-.8cm}\:+\Vert\nabla\bn\Vert_{L^{2p}(D_{9\rho/16}(x))}\big(\Vert\nabla\bR\Vert_{L^{2p}(D_{9\rho/16}(x))}+\Vert\nabla S\Vert_{L^{2p}(D_{9\rho/16}(x))}+\Vert\nabla\bX\Vert_{L^{2p}(D_{9\rho/16}(x))}+\Vert\nabla Y\Vert_{L^{2p}(D_{9\rho/16}(x))}  \big)\nonumber\\[1ex]
&&\hspace{-.8cm}\lesssim\:\rho^{\frac{2}{p}-2}M(M+1)\:.
\end{eqnarray}
As shown in \cite{BR1}, the Gauss map $\bn$ satisfies a perturbed harmonic map equation, namely
\bes
|\Delta\bn|\;\leq\;2\big|\text{div}\big(\bH\wedge\nabla^\perp\bp\big)\big|+\text{O}\big(|\nabla\bn|^2\big)\:.
\ees
Accordingly, from (\ref{nas0}) and (\ref{nadege}), we find
\begin{eqnarray}\label{nadege2}
&&\hspace{-.8cm}\Vert\nabla^2\bn\Vert_{L^{p}(D_{17\rho/32}(x))}\;\;\lesssim\;\;\Vert\text{div}\big(\bH\wedge\nabla^\perp\bp\big)\Vert_{L^p(D_{9\rho/16}(x))}+\Vert\nabla\bn\Vert^2_{L^{2p}(D_{9\rho/16}(x))} +\:\rho^{\frac{2}{p}-2}\Vert\nabla\bn\Vert_{L^{2}(D_{9\rho/16}(x))} \nonumber\\[1ex]
&&\lesssim\:\Vert\text{e}^{\la}\bT\Vert_{L^p(D_{\rho}(x))}+\rho^{\frac{2}{p}-2}(M+1)^2+\rho^{\frac{2}{p}-2}\Vert\nabla\bn\Vert_{L^2(D_\rho(x))}\nonumber\\[1ex]
&&\lesssim\:\rho^{\frac{2}{p}-2}(M+1)^2\:.
\end{eqnarray}
\noindent
From this we deduce
\begin{eqnarray}\label{estimnablan}
\left\{\begin{array}{lclcl}
\Vert\nabla\bn \Vert_{L^{2p/(2-p)}(D_{\rho/2}(x))}&\lesssim&\rho^{\frac{2}{p}-2}(M+1)^2&,&\text{if}\:\:p\in(1,2)\\[1.5ex]
\Vert\nabla\bn \Vert_{L^{q}(D_{\rho/2}(x))}&\lesssim&\rho^{\frac{2}{q}-1}(M+1)^2&,&\forall\:\:q<\infty\:,\:\:\text{if}\:\:p=2\\[1.5ex]
\Vert\nabla\bn \Vert_{L^{\infty}(D_{\rho/2}(x))}&\lesssim& \rho^{-1}(M+1)^2&,&\text{if}\:\:p>2\:.
\end{array}\right.
\end{eqnarray}

To complete the proof of Theorem \ref{estim}-(i), we show that the estimates (\ref{estimH}) and (\ref{estimnablan}) may be slightly improved when $p\ge2$. Firstly, when $p>2$, we see that (\ref{pjo1}) holds for all $s<\infty$. In particular, (\ref{simone}), (\ref{tita0}), and (\ref{nas0}) hold with any $s<\infty$ in place of $2p$. But according to (\ref{estimH}) and (\ref{estimnablan}), we know that $\text{e}^\la\bH$ and $\nabla\bn$ are bounded when $p>2$. We may thus let $s$ tend to infinity to find that
\be\label{keane1}
\Vert\text{e}^{\la}\bH\Vert_{L^\infty(D_{\rho/2}(x))}\;\lesssim\;\rho^{-1}M\qquad\text{and}\qquad \Vert\nabla\bn\Vert_{L^\infty(D_{\rho/2}(x))}\;\lesssim\;\rho^{-1}(M+1)\:.
\ee
Similarly, when $p=2$, we have that (\ref{pjo1}), and thus (\ref{simone}) and (\ref{nas0}) hold for all $s<\infty$ in place of $2p$. In particular,
\bes
\Vert\text{e}^{\la}\bH\Vert_{L^q(D_{\rho/2}(x))}\;\lesssim\;\rho^{\frac{2}{q}-1}M\qquad\text{and}\qquad \Vert\nabla\bn\Vert_{L^q(D_{\rho/2}(x))}\;\lesssim\;\rho^{\frac{2}{q}-1}(M+1)\:,\qquad\forall\:\:q<\infty\:.
\ees


\subsection{Proof of Theorem \ref{estim}-(ii)}

In this section we will build upon the results previously derived in order to obtain regularity estimates for an inhomogeneous Willmore equation of the type
\bes
\text{div}\big(\nabla\bH-2\pro\nabla\bH+|\bH|^2\nabla\bp  \big)\;=\;\bv\:,
\ees
where we suppose that
\bes
\text{e}^{\la}\bv\,\in\,L^r(D_1(0))\qquad\text{for some $r\ge1$}\:.
\ees
Let $D_\rho(x)\subset D_1(0)$. In order to recover (\ref{ourdiv2}), we let $\bV$ satisfy the problem
\be\label{eqv}
\left\{\begin{array}{rclcl}
-\,\Delta\bV&=&\bv\,,&\ &\text{in}\:\:D_\rho(x)\\[1ex]
\bV&=&\vec{0}\,,&\ &\text{on}\:\:\partial D_\rho(x)\:.
\end{array}\right.
\ee
Using the Harnack inequality (\ref{harnacko}), we easily deduce 
\bes
\left\{\begin{array}{lclcl}
\Vert\text{e}^{\la}\nabla\bV\Vert_{L^{2,\infty}(D_\rho(x))}&\lesssim&\Vert\text{e}^{\la}\bv\Vert_{L^1(D_\rho(x))}&,&r=1\\[1.5ex]
\Vert\text{e}^{\la}\nabla\bV\Vert_{L^{2r/(2-r)}(D_\rho(x))}&\lesssim&\Vert\text{e}^{\la}\bv\Vert_{L^r(D_\rho(x))}&,&r\in(1,2)\\[1.5ex]
\Vert\text{e}^{\la}\nabla\bV\Vert_{L^{q}(D_\rho(x))}&\lesssim&\rho^{\frac{2}{q}+1-\frac{2}{r}}\Vert\text{e}^{\la}\bv\Vert_{L^r(D_\rho(x))}&,&\forall\:\:q<\infty\:\:,\:\:r\ge2\:.
\end{array}\right.
\ees
We are now back in the case studied in the previous section with $\bT:=\nabla\bV$. In particular, when $r=1$, we find
\bes
\Vert\text{e}^{\la}\bT\Vert_{L^p(D_\rho(x))}\;\lesssim\;\rho^{\frac{2}{p}-1}\Vert\text{e}^{\la}\bT\Vert_{L^{2,\infty}(D_\rho(x))}\;\lesssim\;\rho^{\frac{2}{p}-1}\Vert\text{e}^{\la}\bv\Vert_{L^1(D_\rho(x))}\qquad\forall\:\:p\in(1,2)\:,
\ees
from which (\ref{estimnablan}) and (\ref{nadege2}) yield
\bes\label{turlululu2}
\Vert\nabla^2\bn\Vert_{L^p(D_{\rho/2}(x))}+\Vert\nabla\bn\Vert_{L^{2p/(2-p)}(D_{\rho/2}(x))}\;\lesssim\;\rho^{\frac{2}{p}-2}(M+1)^2\qquad\forall\:\:p\in(1,2)\:,
\ees
where
\bes
M\;=\;\rho\Vert\text{e}^{\la}\bv\Vert_{L^{1}(D_\rho(x))}+\Vert\nabla\bn\Vert_{L^2(D_\rho(x))}\:.
\ees

On the other hand, when $r>1$, $\text{e}^{\la}\bT$ lies in a space of the type $L^{2+\delta}$ for some $\delta>0$. This time, the estimates (\ref{keane1}) and (\ref{nadege2}) give
\be\label{keane0}
\rho\Vert\text{e}^{\la}\bH\Vert_{L^\infty(D_{\rho/2}(x))}\;\lesssim\;M\:,\qquad\rho\Vert\nabla\bn\Vert_{L^\infty(D_{\rho/2}(x))}\;\lesssim\;1+M\:,
\ee
and
\be\label{turlululu8}
\left\{\begin{array}{rclcl}
\rho^{3-\frac{2}{r}}\Vert\nabla^2\bn\Vert_{L^{2r/(2-r)}(D_{\rho/2}(x))}&\lesssim&(M+1)^2&,&\quad\text{if}\:\:\:r\in(1,2)\\[2.5ex]
\rho^{2-\frac{2}{q}}\Vert\nabla^2\bn\Vert_{L^{q}(D_{\rho/2}(x))}&\lesssim&(M+1)^2&,&\forall\:\:q<\infty\quad\text{if}\:\:\:r\ge2\:,
\end{array}\right.
\ee
with
\be\label{titaminus1}
M\;=\;\rho^{3-\frac{2}{r}}\Vert\text{e}^{\la}\bv\Vert_{L^{r}(D_\rho(x))}+\Vert\nabla\bn\Vert_{L^2(D_\rho(x))}\:.
\ee

We now prove third-order derivative estimates for $\nabla\bn$. For the sake of brevity, we only consider the case $r\in(1,2)$. From (\ref{eqv}), we find
\bes
\Vert\nabla\bV\Vert_{L^{r^*}(D_\rho(x))}+\Vert\nabla^2\bV\Vert_{L^r(D_\rho(x))}\;\lesssim\;\Vert\bv\Vert_{L^r(D_\rho(x))}\:,
\ees
where for notational convenience, we have set $r^*:=2r/(2-r)$. Since $\bT:=\nabla\bV$, the latter and (\ref{harnacko}) yield
\be\label{tita1}
\Vert\text{e}^{\la}\bT\Vert_{L^{r^*}(D_\rho(x))}+\Vert\text{e}^{\la}\nabla\bT\Vert_{L^r(D_\rho(x))}\;\lesssim\;\Vert\text{e}^{\la}\bv\Vert_{L^r(D_\rho(x))}\:.
\ee
On the other hand, using (\ref{tita0}), we have
\be\label{tita2}
\Vert\nabla^2\bp\Vert_{L^{2}(D_{9\rho/16}(x))}\;\lesssim\;\Vert\text{e}^{\la}\Vert_{L^\infty(D_\rho(x))}(1+M)\:,
\ee
where $M$ is as in (\ref{titaminus1}). \\
According to (\ref{eqX}), we find
\be\label{tita4}
\rho\Vert\nabla^2\bX\Vert_{L^{2}(D_{\rho}(x))}+\rho\Vert\nabla^2 Y\Vert_{L^{2}(D_{\rho}(x))}+\Vert\nabla\bX\Vert_{L^{2}(D_{\rho}(x))}+\Vert\nabla Y\Vert_{L^{2}(D_{\rho}(x))}\;\lesssim\;\rho^{2-\frac{2}{r^*}}\Vert\text{e}^{\la}\bv\Vert_{L^r(D_\rho(x))}\:;
\ee
and moreover, owing to (\ref{tita2})
\begin{eqnarray}\label{tita5}
&&\hspace{-2cm}\Vert\nabla\Delta\bX\Vert_{L^r(D_{9\rho/16}(x))}+\Vert\nabla\Delta Y\Vert_{L^r(D_{9\rho/16}(x))}\nonumber\\[1ex]
&&\lesssim\:\:\Vert\nabla^2\bp\Vert_{L^2(D_{9\rho/16}(x))}\Vert\bT\Vert_{L^{r^*}(D_{9\rho/16}(x))}+\Vert\nabla\bp\Vert_{L^\infty(D_{9\rho/16}(x))}\Vert\nabla\bT\Vert_{L^r(D_{9\rho/16}(x))}\nonumber\\[1ex]
&&\lesssim\:\:(1+M)\Vert\text{e}^{\la}\bv\Vert_{L^r(D_\rho(x))}\:.
\end{eqnarray}

Note that (\ref{pjo1}) and (\ref{tita1}) yield
\bes
\Vert\nabla\bR\Vert_{L^2(D_{3\rho/4}(x))}+\Vert\nabla S\Vert_{L^2(D_{3\rho/4}(x))}\;\lesssim\;M\:.
\ees
In addition, (\ref{nadege0}) states
\be\label{tito2}
\Vert\nabla^2\bR\Vert_{L^2(D_{\rho/2}(x))}+\Vert\nabla^2 S\Vert_{L^2(D_{\rho/2}(x))}\;\lesssim\;\rho^{-1}M(M+1)\:.
\ee

From (\ref{backimm}), we easily verify that
\begin{eqnarray*}
\big|\nabla\text{div}\big(\bH\wedge\nabla^\perp\bp\big)\big|&\lesssim&|\nabla\Delta\bX|+|\nabla\Delta Y|+|\nabla\bn|\big(|\nabla^2\bR|+|\nabla^2S|+|\nabla^2\bX|+|\nabla^2Y|  \big)\\
&&\hspace{1cm}+\;\;|\nabla^2\bn|\big(|\nabla\bR|+|\nabla S|+|\nabla\bX|+|\nabla Y|  \big)\:.
\end{eqnarray*}
The estimates (\ref{keane0}), (\ref{turlululu8}), and (\ref{tita4})-(\ref{tito2}) then show that
\begin{eqnarray}
\big\Vert\nabla\text{div}\big(\bH\wedge\nabla^\perp\bp\big)\big\Vert_{L^r(D_{\rho/2}(x))}\;\;\lesssim\;\;\rho^{\frac{2}{r}-3}M(M+1)^2\:.
\end{eqnarray}
We can then proceed as in (\ref{nadege2}) to obtain
\be\label{tito5}
\Vert\nabla^3\bn\Vert_{L^r(D_{\rho/3}(x))}\;\lesssim\;\rho^{\frac{2}{r}-3}(M+1)^3\:.
\ee

The case $r\ge2$ is handled {\it mutatis mutandis}, and one arrives too at the third-order estimate (\ref{tito5}). In particular, when $r>2$, we see that $\nabla^2\bn$ is bounded, and we can slightly improve (\ref{turlululu8}) to
\bes
\rho^2\Vert\nabla^2\bn\Vert_{L^\infty(D_{\rho/3}(x))}\;\lesssim\;(M+1)^2\:.
\ees
This concludes the proof of Theorem \ref{estim}-(ii).

\subsection{On smoothness of the solution: proof of Corollary \ref{smooth}}

Let us suppose that $\bp\in W^{2,2}_{imm}(D_1(0),\R^m)$ satisfies the equation
\bes
\text{div}\big(\nabla\bH-2\pro\nabla\bH+|\bH|^2\nabla\bp+\bT\big)\;=\;\vec{v}\qquad\text{on}\:\:D_1(0)\:,
\ees
where  $\bT$ and $\bv$ are smooth. As we are interested in obtaining a local result, we may always rescale so as to guarantee that the small energy assumption
\bes
\Vert\nabla\bn\Vert_{L^2(D_1(0))}\;<\;\eps_0
\ees
holds for some $\eps_0$ sufficiently small. We proved in the last section that
$\nabla\bn$ is a bounded function. Owing to the Liouville equation\footnote{$K$
denotes the Gauss curvature.}
\bes
-\Delta\la\;=\;\text{e}^{2\la}K\;=\;\text{O}(|\nabla\bn|^2)\:,
\ees
it follows that $\text{e}^{\pm\la}$ lie in $\bigcap_{p<\infty}W^{2,p}$. The function $\bV$ defined in (\ref{eqv}) is smooth. By definition, so is $\bU:=\bT+\nabla\bV$. Using (\ref{eqX}), we deduce that $\bX$ and $Y$ belong to $\bigcap_{p<\infty}W^{4,p}$. We see in the paragraph following (\ref{backimm}) that $\bR$ and $S$ also belong to $\bigcap_{p<\infty}W^{2,p}$. Then the equation (\ref{backimm}) yields now that the immersion lies in $\bigcap_{p<\infty}W^{3,p}$, and thus that $\nabla\bn$ lies in $\bigcap_{p<\infty}W^{2,p}$. This process may now be repeated as much as required to reach the conclusion that $\bp$ is smooth.

\subsection{Remarks about the critical case}\label{remcrit}

As its name indicates, the critical case is far more delicate to handle, and, as far as the authors know, there is no general method to prove the regularity of solutions to the inhomogeneous Willmore equation
\be\label{ourdiv4}
\text{div}\big(\nabla\bH-2\pro\nabla\bH+|\bH|^2\nabla\bp+\bT  \big)\;=\;\vec{0}\:,
\ee
with a generic inhomogeneity $\text{e}^{\la}\bT\in L^1$, if it is only known
that the second fundamental form is square integrable. There are of course
special cases, such as the Willmore immersions (with $\bT\equiv\vec{0}$) and
more generally the conformally constrained Willmore immersions (which include Willmore
and CMC immersions) whose $\text{e}^{\la}\bT$ has a very specific form, see
\cite{Ber1}. The conformally constrained Willmore immersions have an inhomogeneous term $\bT$ for which the solutions to (\ref{eqX}) are identically
vanishing. In turn, this guarantees the system (\ref{sysRS}) is of Wente type
and can thus be made subcritical just as we have done for $\text{e}^{\la}\bT\in
L^{p>1}$. \\

But even if we assume from the onset that the solution to (\ref{ourdiv4}) is
sufficiently regular\footnote{even if ever so slightly, say $\bH\in
L^{2+\delta}$ for some $\delta>0$.}, the presence of an inhomogeneity
$\bT$, and thus of nonzero solutions of (\ref{eqX}), will in general prevent us
from reaching estimates of the type appearing in Theorem \ref{estim}. This
difficulty can only be resolved on a case-by-case basis. We will content
ourselves in this short section with mentioning one specific type of
inhomogeneities for which Theorem \ref{estim} can be obtained. \\
Let us write the inhomogeneity $\bT$ in the form
\bes
\bT\;=\;\left(\begin{array}{c}A_1\\A_2\end{array}\right)\partial_{x^1}\bp+\left(\begin{array}{c}B_1\\B_2\end{array}\right)\partial_{x^2}\bp+\left(\begin{array}{c}\bU_1\\\bU_2\end{array}\right)\;,
\ees
where $\bU_1$ and $\bU_2$ are two normal vectors. One easily verifies that
\bes
\nabla\bp\wedge\bT\;=\;\text{e}^{2\la}(A_2-B_1)(\star\bn)-\bU_1\wedge\partial_{x^1}\bp-\bU_2\wedge\partial_{x^2}\bp
\ees
and
\bes
\nabla\bp\cdot\bT\;=\;\text{e}^{2\la}(A_1+B_2)\:.
\ees
Accordingly, if the functions $(A_1+B_2)$, $(A_2-B_1)$, and the normal projection $\pro\bT$ lie in the space $L^{1+\delta}$ for some $\delta>0$, we can apply to (\ref{eqX}) the same technique as that used in the proof of Theorem \ref{estim}. This holds of course even if the functions $A_1$, $A_2$, $B_1$, and $B_2$ are only merely integrable.  \\

In general, it is not possible to obtain a subcritical-type energy estimate. However, as we have seen above, there are exceptions when $\bT$ has a specific form. Another important exceptional case occurs when $\bT$ depends on the geometry of the problem, and if the solution is already known to be regular enough, say $\bp\in W^{2,2+\delta}(D_1(0))$, for some positive $\delta\in(0,1)$. We only focus on the specific situation when
\bes
\text{e}^{\la}\bT\;=\;\text{O}\big(|\nabla\bn|^2\big)
\ees
for an inhomogeneous Willmore problem of the type
\bes
\text{div}\big(\nabla\bH-2\pro\nabla\bH+|\bH|^2\nabla\bp+\bT  \big)\;=\;\vec{0}\:,
\ees
and assuming as usual that
\bes
\Vert\nabla\bn\Vert_{L^2(D_1(0))}\;<\;\eps_0
\ees
for some $\eps_0$ chosen sufficiently small.\\ 
As $\text{e}^{\la}\bT\in L^{1+\frac{\delta}{2}}$, it follows from Theorem \ref{estim}, that $\nabla\bn$ lies in $W^{1,1+\frac{\delta}{2}}\subset L^{2\frac{2+\delta}{2-\delta}}$, and thus $\text{e}^{\la}\bT$ lies in $W^{1,2\frac{2+\delta}{6-\delta}}$, which is a proper subset of $L^{1+\frac{\delta}{2}}$. Calling again on Theorem \ref{estim}, the integrability of $\nabla\bn$ is improved accordingly. This procedure may be repeated until reaching that $\nabla^2\bn$ belongs to all $L^{p}$ spaces, with $p$ finite, i.e. that $\bn$ belongs to $C^{1,\al}$ for all $\al<1$. Standard arguments then imply that $\bn$, and thus the immersion $\bp$, are smooth.


\subsection{Gap phenomenon: proof of Theorem \ref{gap0}}

Let us suppose that $\Sigma$ is a complete, connected, non-compact, oriented,  immersed surface into $\mathbb{R}^{m\ge3}$ satisfying an inhomogeneous Willmore equation (\ref{willeq}) of the form 
\be\label{rozen5}
\Delta_\perp\bH+\big\langle\bA\cdot\bH,\bA\big\rangle_g-2|\bH|^2\bH\;=\;\vec{\mathcal{W}}\:,
\ee
with the same notation as before, and where $\vec{\mathcal{W}}$ is a normal field with the property that
\be\label{rozen6}
\vec{\mathcal{W}}\;=\;\text{O}(|\bA|^3)\qquad\text{i.e.}\qquad \dfrac{1}{c}|\bA|^3\leq|\vec{\mathcal{W}}|\leq c|\bA|^3\:,
\ee
for some constant $c\ge1$. We suppose further that 
\be\label{minus0}
\int_\Sigma|\bA|_g^2\,d\text{vol}_g\;<\;\eps_0^2\:,
\ee
for some $\eps_0^2$ chosen to be small enough (at least smaller than $8\pi/3$). A well-known result of M\"uller and Sverak \cite{MS} guarantees that $\Sigma$ is embedded and conformally equivalent to $\R^2$. Accordingly, we parametrize $\Sigma$ by a conformal immersion $\bp:\R^2\hookrightarrow\R^{m}$ with conformal parameter $\la$, and such that $\bp\in W^{2,2}(\R^2)$. \\
Just as was done in Section \ref{intro}, in the flat coordinates of $\R^2$, the inhomogeneous Willmore equation (\ref{rozen5}) can be recast in the form 
\bes
\text{div}\big(\nabla\bH-2\pro\nabla\bH+|\bH|^2\nabla\bp   \big)\;=\;\bv\qquad\text{on}\:\:\R^2\:,
\ees
where
\bes
\bv\;:=\;\text{e}^{2\la}\vec{\mathcal{W}}\:.
\ees
Per (\ref{rozen6}), note that
\be\label{rozen7}
|\text{e}^{\la}\bv|\;\simeq\;|\text{e}^{\la}\bA|^3\;\simeq\;|\nabla\bn|^3\:,
\ee
where, as always, $\bn$ is the Gauss map associated with $\bp$. The smallness hypothesis (\ref{minus0}) translates into
\be\label{minus}
\Vert\nabla\bn\Vert_{L^2(\R^2)}\;<\;\eps_0\:,
\ee
for some $\eps_0>0$ sufficiently small.\\

Owing to the Liouville equation
\bes
-\Delta\la\;=\;\text{e}^{2\la}K\;=\;\text{O}(|\nabla\bn|^2)\,\in\,L^1(\R^2)\:,
\ees
it follows that $\nabla\la$ lies in the space $L^{2,\infty}(\R^2)$ with norm controlled by $\Vert\nabla\bn\Vert_{L^2(\R^2)}$. We can in particular repeat the analysis leading to Lemma \ref{harnack} to deduce that
\be\label{zuzu0}
\Vert\text{e}^{\la}\Vert_{L^\infty(\R^2)}\Vert\text{e}^{-\la}\Vert_{L^\infty(\R^2)}\;\leq\;C(\eps_0)\:.
\ee
This Harnack-type inequality will be used in our argument. \\

As we did in the proof of Theorem \ref{estim}-(ii), we let
\bes
-\Delta\bV\;=\;\bv\qquad\text{on}\:\:\R^2\:,
\ees
and $\bT:=\nabla\bV$.
As the equation
\bes
\text{div}\big(\nabla\bH-2\pro\nabla\bH+|\bH|^2\nabla\bp+\bT   \big)\;=\;\vec{0}\qquad\text{holds on}\:\:\R^2\:,
\ees
we can repeat the analysis done in the proof of Theorem \ref{estim}-(i) and deduce the existence of $\bR$ and $S$ satisfying 
\be\label{sysRSquad}
\left\{\begin{array}{lcl}
\Delta\bR&=&\nabla(\star\bn)\bullet\nabla^\perp\bR+\nabla(\star\bn)\cdot\nabla^\perp S+\text{div}\big((\star\bn)\bullet\nabla\bX+(\star\bn)\nabla Y   \big)\\[1ex]
\Delta S&=&\nabla(\star\bn)\cdot\nabla^\perp\bR+\text{div}\big((\star\bn)\cdot\nabla\bX  \big)\:,
\end{array}\right.
\ee
where, as before, $\bX$ and $Y$ satisfy 
\be\label{eqXquad}
\Delta\bX\;=\;\nabla\bp\wedge\bT\qquad\text{and}\qquad\Delta Y\;=\;\nabla\bp\cdot\bT\qquad\text{on}\:\:\R^2\:.
\ee
We have
\be\label{ziba}
\Vert\Delta\bX\Vert_{L^q(\R^2)}+\Vert\Delta Y\Vert_{L^q(\R^2)}+\Vert\nabla\bX\Vert_{L^{q^*}(\R^2)}+\Vert\nabla Y\Vert_{L^{q^*}(\R^2)}\;\lesssim\;\Vert\text{e}^{\la}\bT\Vert_{L^{q}(\R^2)}\:,
\ee
for $q\in(1,2)$ and $q^*:=2q/(2-q)$. \\
Applying Wente's inequality to (\ref{sysRSquad}) as in Lemma IV.2 of \cite{BR2}, we find
\bes
\Vert\nabla\bR\Vert_{L^{q^*}(\R^2)}+\Vert\nabla S\Vert_{L^{q^*}(\R^2)}\;\leq\;\Vert\nabla\bn\Vert_{L^2(\R^2)}\big(\Vert\nabla\bR\Vert_{L^{q^*}(\R^2)}+\Vert\nabla S\Vert_{L^{q^*}(\R^2)}\big)+\Vert\nabla\bX\Vert_{L^{q^*}(\R^2)}+\Vert\nabla Y\Vert_{L^{q^*}(\R^2)}\:,
\ees
which, owing to (\ref{minus}) and (\ref{ziba}), yields
\be\label{ziba2}
\Vert\nabla\bR\Vert_{L^{q^*}(\R^2)}+\Vert\nabla S\Vert_{L^{q^*}(\R^2)}\;\leq\;C(\eps_0)\big(\Vert\nabla\bX\Vert_{L^{q^*}(\R^2)}+\Vert\nabla Y\Vert_{L^{q^*}(\R^2)}\big)\;\leq\;C(\eps_0)\Vert\text{e}^{\la}\bT\Vert_{L^{q}(\R^2)}\:.
\ee
We have seen in the previous section that
\bes
2\bH\wedge\nabla^\perp\bp\;=\;\big(\nabla^\perp\bR-\nabla\bX)\bul(\star\bn)+\big(\nabla^\perp S-\nabla Y)(\star\bn)\:,
\ees
hence
\begin{eqnarray*}
\big|\text{div}(2\bH\wedge\nabla^\perp\bp)\big|&\leq&|\Delta\bX|+|\Delta Y|+|\nabla\bn|\big(|\nabla\bR|+|\nabla S|+|\nabla\bX|+|\nabla Y|\big)\:.
\end{eqnarray*}
This gives, using (\ref{ziba}) and (\ref{ziba2}),
\begin{eqnarray}\label{rozen3}
\big\Vert\text{div}(2\bH\wedge\nabla^\perp\bp)\big\Vert_{L^q(\R^2)}&\leq&\Vert\Delta\bX\Vert_{L^q(\R^2)}+\Vert\Delta Y\Vert_{L^q(\R^2)}\nonumber\\
&&\hspace{-2cm}+\:\Vert\nabla\bn\Vert_{L^2(\R^2)}\big(\Vert\nabla\bR\Vert_{L^{q^*}(\R^2)}+\Vert\nabla S\Vert_{L^{q^*}(\R^2)}+\Vert\nabla\bX\Vert_{L^{q^*}(\R^2)}+\Vert\nabla Y\Vert_{L^{q^*}(\R^2)}\big)\nonumber\\[1ex]
&&\hspace{-4cm}\leq\:2\Vert\text{e}^{\la}\bT\Vert_{L^q(\R^2)}+C(\eps_0)\big( \Vert\nabla\bX\Vert_{L^{q^*}(\R^2)}+\Vert\nabla Y\Vert_{L^{q^*}(\R^2)} \big)\nonumber\\[1ex]
&&\hspace{-4cm}\leq\:C(\eps_0)\Vert\text{e}^{\la}\bT\Vert_{L^q(\R^2)}\:.
\end{eqnarray}
Recall that
\bes
\Delta\bn\;=\;\text{div}\big(2\bH\wedge\nabla^\perp\bp\big)+\text{O}(|\nabla\bn|^2)\:.
\ees
According to (\ref{rozen3}), to (\ref{minus}), and to the Sobolev embedding theorem, we thus have
\begin{eqnarray*}
\Vert\nabla^2\bn\Vert_{L^q(\R^2)}&\leq&\big\Vert\text{div}(2\bH\wedge\nabla^\perp\bp)\big\Vert_{L^q(\R^2)}+\Vert\nabla\bn\Vert_{L^2(\R^2)}\Vert\nabla\bn\Vert_{L^{q^*}(\R^2)}\nonumber\\[1ex]
&\leq&C(\eps_0)\Vert\text{e}^{\la}\bT\Vert_{L^q(\R^2)}+\eps_0\Vert\nabla^2\bn\Vert_{L^q(\R^2)}\:,
\end{eqnarray*}
thereby yielding
\be\label{rozen4}
\Vert\nabla^2\bn\Vert_{L^q(\R^2)}\;\leq\;C(\eps_0)\Vert\text{e}^{\la}\bT\Vert_{L^q(\R^2)}\:.
\ee
We now call upon the Gagliardo-Nirenberg interpolation inequality, (\ref{minus}), and (\ref{rozen4})  to find
\bes
\Vert\nabla\bn\Vert_{L^{p}(\R^2)}\;\leq\;\Vert\nabla^2\bn\Vert^{\al}_{L^q(\R^2)}\Vert\nabla\bn\Vert^{1-\al}_{L^2(\R^2)}\;\leq\;C(\eps_0)\eps_0^{1-\al}\Vert\text{e}^{\la}\bT\Vert^{\al}_{L^q(\R^2)}\:,
\ees
for
\bes
\dfrac{1}{p}\;=\;\dfrac{1}{2}+\bigg(\dfrac{1}{q}-1\bigg)\al\qquad\text{and}\qquad 0\leq\al\leq1\:.
\ees
Equivalently, 
\bes
\big\Vert|\nabla\bn|^3\big\Vert_{L^{b}(\R^2)}\;\leq\;C(\eps_0)\eps_0^{3(1-\al)}\Vert\text{e}^{\la}\bT\Vert^{3\al}_{L^q(\R^2)}\:,
\ees
for
\bes
\dfrac{1}{b}\;=\;\dfrac{3}{2}+3\bigg(\dfrac{1}{q}-1\bigg)\al\:.
\ees
As $\text{e}^\la\Delta\bV=-\text{e}^{\la}\bv=\text{O}(|\nabla\bn|^3)$ and $\bT=\nabla\bV$, the latter yields
\bes
\Vert\text{e}^{\la}\Delta\bV\Vert_{L^{b}(\R^2)}\;\leq\;C(\eps_0)\eps_0^{3(1-\al)}\Vert\text{e}^{\la}\nabla\bV\Vert^{3\al}_{L^q(\R^2)}\:,
\ees
hence, using (\ref{zuzu0}), 
\be\label{zuzu1}
\Vert\Delta\bV\Vert_{L^{b}(\R^2)}\;\leq\;C(\eps_0)\eps_0^{3(1-\al)}\Vert\text{e}^{\la}\Vert^{3\al-1}_{L^\infty(\R^2)}\Vert\nabla\bV\Vert^{3\al}_{L^q(\R^2)}\:.
\ee
Let $\delta\in(0,2/3)$. We specialize to
\bes
q\;=\;2-\delta\qquad\text{and}\qquad 3\al\;=\;\dfrac{1}{1-\delta}\:.
\ees
This gives
\bes
\dfrac{1}{b}\;=\;\dfrac{3}{2}-\dfrac{1}{2-\delta}\:,\qquad\text{so that}\qquad b\in(1,2)\:.
\ees
Using the Sobolev embedding theorem in (\ref{zuzu1}) then gives
\be\label{zuzu2}
\Vert\nabla\bV\Vert^{1-\delta}_{L^{\frac{2-\delta}{1-\delta}}(\R^2)}\;\leq\;C(\eps_0)\eps_0^{{2-3\delta}}\Vert\text{e}^{\la}\Vert^{{\delta}}_{L^\infty(\R^2)}\Vert\nabla\bV\Vert_{L^{2-\delta}(\R^2)}\:.
\ee
Since 
\bes
\dfrac{1-\delta}{2-\delta}+\dfrac{1}{2-\delta}\;=\;1\:,
\ees
we interpolate (\ref{zuzu2}) to find
\bes
\Vert\nabla\bV\Vert^{2(1-\delta)}_{L^2(\R^2)}\;\leq\;C(\eps_0)\eps_0^{{2-3\delta}}\Vert\text{e}^{\la}\Vert^{{\delta}}_{L^\infty(\R^2)}\Vert\nabla\bV\Vert^{2-\delta}_{L^{2-\delta}(\R^2)}\:.
\ees
Letting $\delta\searrow0$ reveals that
\bes
\Vert\nabla\bV\Vert_{L^2(\R^2)}\;\leq\;C(\eps_0)\eps_0\Vert\nabla\bV\Vert_{L^{2}(\R^2)}\:.
\ees
Since $\eps_0$ can be adjusted at will, the latter implies that $\nabla\bV\equiv\vec{0}$, hence that $\bv=-\Delta\bV\equiv\vec{0}$, and therefore that $\nabla\bn=\vec{0}$. This guarantees that $\Sigma$ is a flat plane, as announced.

\bigskip

\renewcommand{\theequation}{A.\arabic{equation}}
\renewcommand{\theTh}{A.\arabic{Th}}
\renewcommand{\theProp}{A.\arabic{Prop}}
\renewcommand{\theLma}{A.\arabic{Lma}}
\renewcommand{\theLm}{A.\arabic{Lm}}
\renewcommand{\theCo}{A.\arabic{Co}}
\renewcommand{\theRm}{A.\arabic{Rm}}
\renewcommand{\theequation}{A.\arabic{equation}}
\setcounter{equation}{0} 
\reset
\appendix
\section{Appendix}
\subsection{Notational conventions}

We append an arrow to all the elements belonging to $\R^m$. To simplify the notation, by $\bp\in X(\di)$ is meant $\bp\in X(\di,\R^m)$ whenever $X$ is a function space. Similarly, we write $\nabla\bp\in X(\di)$ for $\nabla\bp\in \mathbb{R}^2\otimes X(\di,\R^m)$.\\[1.5ex]
We let differential operators act on elements of $\R^m$ componentwise. Thus, for example, $\nabla\bp$ is the element of $\R^2\otimes\R^m$ with $\R^m$-valued components $(\partial_{x^1}\bp,\partial_{x^2}\bp)$. If $S$ is a scalar and $\bR$ an element of $\R^m$, then we let
\begin{eqnarray*}
\bR\cdot\nabla\bp&:=&\big(\bR\cdot\partial_{x^1}\bp\,,\,\bR\cdot\partial_{x^2}\bp\big)\:\\[1ex]
\nabla^\perp S\cdot\nabla\bp&:=&\partial_{x^1} S\,\partial_{x^2}\bp\,-\,\partial_{x^2} S\,\partial_{x^1}\bp\:\\[1ex]
\nabla^\perp\bR\cdot\nabla\bp&:=&\partial_{x^1}\bR\cdot\partial_{x^2}\bp\,-\,\partial_{x^2}\bR\cdot\partial_{x^1}\bp\:\\[1ex]
\nabla^\perp\bR\wedge\nabla\bp&:=&\partial_{x^1}\bR\wedge\partial_{x^2}\bp\,-\,\partial_{x^2}\bR\wedge\partial_{x^1}\bp\:.
\end{eqnarray*}
Analogous quantities are defined according to the same logic. \\

Two operations between multivectors are useful. The {\it interior multiplication} $\res$ maps a pair comprising a $q$-vector $\gamma$ and a $p$-vector $\beta$ to a $(q-p)$-vector. It is defined via
\bes
\langle \gamma\res\beta\,,\alpha\rangle\;=\;\langle \gamma\,,\beta\wedge\alpha\rangle\:\qquad\text{for each $(q-p)$-vector $\al$.}
\ees
Let $\al$ be a $k$-vector. The {\it first-order contraction} operation $\bul$ is defined inductively through 
\bes
\al\bul\beta\;=\;\al\res\beta\:\:\qquad\text{when $\beta$ is a 1-vector}\:,
\ees
and
\bes
\al\bul(\beta\wedge\gamma)\;=\;(\al\bul\beta)\wedge\gamma\,+\,(-1)^{pq}\,(\al\bul\gamma)\wedge\beta\:,
\ees
when $\beta$ and $\gamma$ are respectively a $p$-vector and a $q$-vector. 

\subsection{Some useful elliptic results}

The following result is established in the Appendix of \cite{BR0}. 

\begin{Lm}\label{L2weak}
Let $D$ be a disk and suppose that $G=G_1+G_2$  satisfies
\bes
div\,G\;=\;0\qquad\text{on}\quad D\:,
\ees
where
\bes
G_1\,\in\,W^{-1,2}(D,\mathbb{R}^2)\quad,\quad G_3\,\in\,L^1(D,\mathbb{R}^2)\:.
\ees
Then there exists an element $L$ in the space $L^{2,\infty}(D,\mathbb{R})$ such that
\bes
G\;=\;\nabla^{\perp}L\:,
\ees
and
\bes
\Vert L-L_D\Vert_{L^{2,\infty}(D)}\;\leq\;C\big(\Vert G_1\Vert_{W^{-1,2}(D)}+\Vert G_2\Vert_{L^1(D)}\big)\:,
\ees
where $L_D$ denotes the average of $L$ on the disk $D$, and $C$ is a universal constant. 
\end{Lm}

\medskip

\begin{Prop}\label{morreydecay}
Let $D_\rho(x)\subset D_1(0)$, and let $u\in W^{1,(2,\infty)}(D_\rho(x))$ satisfy the equation
\be\label{equ}
\Delta u\:=\:\nabla b\cdot\nabla^\perp u\,+\text{div}\,(b\,\nabla f)\qquad\quad\text{on}\:\:\:\:D_\rho(x)\:,
\ee
where $\,f\in W_0^{2,p}(D_\rho(x))$ for some $p>1$, with
\bes
\Vert\nabla^2 f\Vert_{L^p(D_\rho(x))}\;\lesssim\;\Vert F\Vert_{L^p(D_\rho(x))}\:.
\ees
Suppose moreover that
\be\label{hypn}
b\,\in\,W^{1,2}\cap L^{\infty}(D_\rho(x))\qquad\text{with}\qquad\Vert\nabla b\Vert_{L^2(D_\rho(x))}\;<\;\eps_0\quad\text{and}\quad \Vert b\Vert_{L^\infty(D_\rho(x))}\;\leq\;1\:,
\ee
for some $\eps_0$ chosen to be ``small enough". Then
\bes
\Vert\nabla u\Vert_{L^{s}(D_{5\rho/8}(x))}\;\leq\;C(\eps_0)\Big[\rho^{\frac{2}{s}-1}\Vert\nabla u\Vert_{L^{2,\infty}(D_{\rho}(x))}+\rho^{\frac{2}{s}-\frac{2}{p}+1}\Vert F\Vert_{L^{p}(D_{\rho}(x))}\Big]\:,
\ees
for some constant $C(\eps_0)$ depending only on $\eps_0$, and where $s<2/(2-p)$ if $p\in(1,2)$, or $s<\infty$ if $p\ge2$. 
\end{Prop}
$\textbf{Proof.}$ Suppose first that $p\in(1,2)$. Then for every $D_\sigma(z)\subset D_\rho(x)$, it holds
\be\label{fried}
\Vert\nabla f\Vert_{L^2(D_\sigma(z))}\;\lesssim\;\sigma^{2-\frac{2}{p}}\Vert\nabla f\Vert_{L^{2p/(2-p)}(D_\rho(x))}\;\lesssim\;\sigma^{2-\frac{2}{p}}\Vert F\Vert_{L^p(D_\rho(x))}\:.
\ee
Let us fix once and for all some point $x_0\in D_{3\rho/4}(x)$ and some radius $\,0<r\leq\rho/4$, so that the disk $D_{r}(x_0)$ of radius $r$ and centered on the point $x_0$ is contained in $D_\rho(x)$. With the help of the theorem of Fubini, we may always find some $r_0\in(r/2\,,r)$ such that
\begin{eqnarray}\label{solersestuncon}
\int_{\partial D_{r_0}(x_0)}|\nabla u|^{\frac{3}{2}}&\lesssim&\dfrac{1}{r}\,\int_{D_r(x_0)}|\nabla u|^{\frac{3}{2}}\;\;\lesssim\;\; r^{-\frac{1}{2}}\,\Vert\nabla u\Vert^{\frac{3}{2}}_{L^{2,\infty}(D_r(x_0))}\;\;\lesssim\;\; r_0^{-\frac{1}{2}}\,\Vert\nabla u\Vert^{\frac{3}{2}}_{L^{2,\infty}(D_\rho(x))}\:.
\end{eqnarray}
We next define $u=u_0+u_1$, where the new variables, in accordance with (\ref{equ}), satisfy
\bes
\left\{\begin{array}{rcl}
\Delta u_0&=&\text{div}\,(b\,\nabla f)\\[1ex]
u_0&=&u\end{array}\right.\quad,\quad
\left.\begin{array}{rclcl}
\Delta u_1&=&\nabla b\cdot\nabla^\perp u_{}&\quad&\text{in\:\:\:} D_{r_0}(x_0)\\[1ex]
u_{1}&=&0&\quad&\text{on\:\:\:}\partial D_{r_0}(x_0)\:.
\end{array}\right.
\ees
Let
\bes
\bar{u}\;:=\;\dfrac{1}{2\pi r_0}\int_{\partial D_{r_0}(x_0)}u\:.
\ees\\
Standard elliptic theory, our assumptions on $b$ and $f$, and the Sobolev embedding theorem give
\begin{eqnarray}\label{atchoum}
\Vert\nabla u_0\Vert_{L^2(D_{r_0}(x_0))}&\lesssim&\Vert b\nabla f\Vert_{L^2(D_{r_0}(x_0))}+\Vert u-\bar{u}\Vert_{H^{1/2}(\partial D_{r_0}(x_0))}\nonumber\\[1ex]
&\lesssim&\Vert\nabla f\Vert_{L^2(D_{r_0}(x_0))}+r_0^{\frac{1}{3}}\Vert \nabla u\Vert_{L^{3/2}(\partial D_{r_0}(x_0))}\nonumber\\[1ex]
&\lesssim&r_0^{2-\frac{2}{p}}\Vert F\Vert_{L^p(D_\rho(x))}+\Vert \nabla u\Vert_{L^{2,\infty}(D_\rho(x))}\:,
\end{eqnarray}
where (\ref{fried}) and (\ref{solersestuncon}) were used. \\
To handle $u_1$, we apply Wente's inequality in the form of Lemma IV.2 in \cite{BR2} to obtain
\begin{eqnarray}\label{b20}
\Vert\nabla u_1\Vert_{L^2(D_{r_0}(x_0))}&\lesssim&\Vert\nabla b\Vert_{L^2(D_{r_0}(x_0))}\,\Vert\nabla u\Vert_{L^{2,\infty}(D_{r_0}(x_0))}\;\;\leq\;\;\eps_0\,\Vert\nabla u\Vert_{L^{2,\infty}(D_\rho(x))}\:.
\end{eqnarray}
Altogether, (\ref{atchoum}) and (\ref{b20}) yield that $\nabla u$ belongs to $L^2(D_{r_0}(x_0))$. In particular
\be\label{pj77}
\Vert\nabla u\Vert_{L^2(D_{r_0}(x_0))}\;\lesssim\;r_0^{2-\frac{2}{p}}\Vert F\Vert_{L^p(D_\rho(x))}+\Vert \nabla u\Vert_{L^{2,\infty}(D_\rho(x))}\:.
\ee

Let now $k\in(0,1)$. Using again (\ref{fried}) and standard elliptic theory and growth estimates give
\begin{eqnarray}\label{newkid}
\Vert\nabla u_0\Vert_{L^2(D_{kr_0}(x_0))}&\lesssim&\Vert b\nabla f\Vert_{L^2(D_{r_0}(x_0))}+k\Vert \nabla u_0\Vert_{L^2(D_{r_0}(x_0))}\nonumber\\[1ex]
&\lesssim&\Vert \nabla f\Vert_{L^2(D_{r_0}(x_0))}+k\Vert \nabla u_0\Vert_{L^2(D_{r_0}(x_0))}\nonumber\\[1ex]
&\lesssim&r_0^{2-\frac{2}{p}}\Vert F\Vert_{L^p(D_{\rho}(x))}+k\Vert \nabla u_0\Vert_{L^2(D_{r_0}(x_0))}\:.
\end{eqnarray}
For $u_1$, we apply Wente's inequality this time as in Theorem 3.4.1 of \cite{Hel} so as to find
\be\label{b2}
\Vert\nabla u_1\Vert_{L^2(D_{r_0}(x_0))}\;\lesssim\;\Vert\nabla b\Vert_{L^2(D_{r_0}(x_0))}\,\Vert\nabla u_{}\Vert_{L^2(D_{r_0}(x_0))}\;\lesssim\;\eps_0\,\Vert\nabla u_{}\Vert_{L^2(D_{r_0}(x_0))}\:,
\ee
again up to some multiplicative constant without bearing on the sequel. Hence, combining (\ref{newkid}) and (\ref{b2}) we obtain the estimate
\begin{eqnarray}\label{prechi}
\Vert\nabla u\Vert_{L^2(D_{kr_0}(x_0))}&\lesssim&(k+\eps_0+k\eps_0)\,\Vert\nabla u\Vert_{L^2(D_{r_0}(x_0))}\,+\,r_0^{2-\frac{2}{p}}\Vert F\Vert_{L^{p}(D_{\rho}(x))}\:.
\end{eqnarray}
Because $\eps_0$ is a small adjustable parameter, we may always choose $k$ so as to arrange for $(k+\eps_0+k\eps_0)$ to be small enough. A standard controlled-growth argument (see Lemma III.2.1 in \cite{Gia}) along with (\ref{pj77}) enables us to conclude that for some constant $C(\eps_0)$, there holds
\begin{eqnarray*}
\Vert\nabla u\Vert_{L^2(D_{\sigma}(x_0))}&\leq&C(\eps_0)\sigma^{2-\frac{2}{p}}\Big[r_0^{\frac{2}{p}-2}\Vert\nabla u\Vert_{L^2(D_{r_0}(x_0))}\,+\,\Vert F\Vert_{L^{p}(D_{\rho}(x))}\Big]\nonumber\\[1ex]
&\leq&C(\eps_0)\sigma^{2-\frac{2}{p}}\Big[r_0^{\frac{2}{p}-2}\Vert\nabla u\Vert_{L^{2,\infty}(D_{\rho}(x))}\,+\,\Vert F\Vert_{L^{p}(D_{\rho}(x))}\Big]\:,
\end{eqnarray*}
for
\bes
x_0\in D_{3\rho/4}(x)\qquad\text{and}\qquad \sigma\in(0,r_0)\:.
\ees
In particular, for $r_0=\rho/4$, we find
\be\label{decu}
\Vert\nabla u\Vert_{L^2(D_{\sigma}(x_0))}\;\leq\;C(\eps_0)\sigma^{2-\frac{2}{p}}\Big[\rho^{\frac{2}{p}-2}\Vert\nabla u\Vert_{L^{2,\infty}(D_{\rho}(x))}\,+\,\Vert F\Vert_{L^{p}(D_{\rho}(x))}\Big]\:,
\ee
for
\bes
x_0\in D_{3\rho/4}(x)\qquad\text{and}\qquad \sigma\in(0,\rho/4)\:.
\ees

\medskip
Consider next the maximal function
\be\label{maxfun}
Mg(y)\;:=\;\sup_{\sigma>0}\;\sigma^{\frac{2}{p}-2}\!\int_{D_{\sigma}(y)}|g(z)|\,dz\:.
\ee
We recast the equation (\ref{equ}) in the form
\bes
-\,\Delta u\;=\;b\,\Delta f\;+\;\nabla b\cdot\big(\nabla^\perp u+\nabla f\big)\:.
\ees
Calling upon (\ref{hypn})-(\ref{fried}) and upon the estimate (\ref{decu}), we derive that for $y\in D_{3\rho/4}(x)$, there holds
\begin{eqnarray}\label{prechy}
&&\hspace{-1cm}M\big(\chi_{D_{\rho/4}(y)}\Delta u(z)\big)(y)\:\:\leq\:\:\Vert b\Vert_{L^\infty(D_{\rho}(x))}\,\sup_{0<\sigma<\frac{\rho}{4}}\,\sigma^{\frac{2}{p}-2}\Vert\Delta f\Vert_{L^1(D_{\sigma}(y))}\nonumber\\
&&\hspace{3cm}+\:\:C(\eps_0)\Vert\nabla b\Vert_{L^2(D_{\rho}(x))}\,\sup_{0<\sigma<\frac{\rho}{4}}\sigma^{\frac{2}{p}-2}\Big(\Vert\nabla u\Vert_{L^2(D_{\sigma}(y))}\,+\,\Vert\nabla f\Vert_{L^2(D_{\sigma}(y))}\Big)  \nonumber\\[1ex]
&&\hspace{2.15cm}\leq\:\:C(\eps_0)\Big[\Vert\nabla^2 f\Vert_{L^p(D_{\rho}(x))}+\rho^{\frac{2}{p}-2}\Big(\Vert\nabla u\Vert_{L^{2,\infty}(D_{\rho}(x))}\,+\,\Vert\nabla f\Vert_{L^2(D_{\rho}(x))}\Big)\Big]  \nonumber\\[1ex]
&&\hspace{2.15cm}\leq\:\:C(\eps_0)\Big[\Vert F\Vert_{L^p(D_\rho(x))}\;+\;\rho^{\frac{2}{p}-2}\Vert\nabla u\Vert_{L^{2,\infty}(D_{\rho}(x))}\Big]\:.
\end{eqnarray}
On the other hand, from (\ref{fried}) and (\ref{decu}), we have
\begin{eqnarray}\label{precha}
\Vert\Delta u\Vert_{L^1(D_{3\rho/4}(x))}&\lesssim&\Vert\Delta f\Vert_{L^1(D_{3\rho/4}(x))}+\Vert\nabla u\Vert_{L^2(D_{3\rho/4}(x))}+\Vert\nabla f\Vert_{L^2(D_{3\rho/4}(x))}\nonumber\\[1ex]
&\leq&C(\eps_0)\Big[\rho^{2-\frac{2}{p}}\Vert F\Vert_{L^p(D_{\rho}(x))}+\Vert\nabla u\Vert_{L^{2,\infty}(D_{\rho}(x))}\Big]\:.
\end{eqnarray}
Proposition 3.2 from \cite{Ad} states that
\bes
\big\Vert|z|^{-1}*\chi_{D_{\rho/4}(y)}\Delta u\big\Vert_{L^{\al,\infty}(D_{3\rho/4}(x))}\;\lesssim\;\big\Vert M(\chi_{D_{\rho/4}(y)}\Delta u)\big\Vert^{1-\frac{1}{\al}}_{L^{\infty}(D_{3\rho/4}(x))}\big\Vert\Delta u\big\Vert^{\frac{1}{\al}}_{L^{1}(D_{3\rho/4}(x))}\:,
\ees
where $\al:=2/(2-p)>2$. Hence, according to (\ref{prechy}) and (\ref{precha}), we have
\begin{eqnarray}\label{precho}
&&\hspace{-1.5cm}\big\Vert|z|^{-1}*\chi_{D_{\rho/4}(y)}\Delta u\big\Vert_{L^{\al,\infty}(D_{3\rho/4}(x))}\nonumber\\[1ex]
&&\hspace{-0.75cm}\leq\:C(\eps_0)\Big[\Vert F\Vert_{L^p(D_\rho(x))}\;+\;\rho^{\frac{2}{p}-2}\Vert\nabla u\Vert_{L^{2,\infty}(D_{\rho}(x))}\Big]^{1-\frac{1}{\al}}\Big[\rho^{2-\frac{2}{p}}\Vert F\Vert_{L^p(D_{\rho}(x))}+\Vert\nabla u\Vert_{L^{2,\infty}(D_{\rho}(x))}\Big]^{\frac{1}{\al}}\nonumber\\[1ex]
&&\hspace{-0.75cm}\leq\:C(\eps_0)\rho^{(\frac{2}{p}-2)(1-\frac{1}{\al})}\Big[\rho^{2-\frac{2}{p}}\Vert F\Vert_{L^p(D_\rho(x))}\;+\;\Vert\nabla u\Vert_{L^{2,\infty}(D_{\rho}(x))}\Big]\nonumber\\[1ex]
&&\hspace{-0.75cm}\leq\:C(\eps_0)\rho^{1-p}\Big[\rho^{2-\frac{2}{p}}\Vert F\Vert_{L^p(D_\rho(x))}\;+\;\Vert\nabla u\Vert_{L^{2,\infty}(D_{\rho}(x))}\Big]\:.
\end{eqnarray}

\medskip
We let $y\in D_{3\rho/4}(x)$ and we again decompose $u=u_2+u_3$ with 
\bes
\left\{\begin{array}{rcl}
\Delta u_2&=&0\\[1ex]
u_2&=&u\end{array}\right.\quad,\quad
\left.\begin{array}{rclcl}
\Delta u_3&=&\chi_{D_{\rho/4}(y)}\Delta u&\quad&\text{in\:\:\:} D_{3\rho/4}(x)\\[1ex]
u_{3}&=&0&\quad&\text{on\:\:\:}\partial D_{3\rho/4}(x)\:.
\end{array}\right.
\ees
Let $s\in(2,\al)$. Using standard estimates for the harmonic function $u_2$ and the estimate (\ref{precho}) gives
\begin{eqnarray}\label{bell1}
\Vert\nabla u\Vert_{L^s(D_{5\rho/8}(x))}&\leq&\Vert\nabla u_2\Vert_{L^{s}(D_{5\rho/8}(x))}+\Vert\nabla u_3\Vert_{L^{s}(D_{5\rho/8}(x))}\nonumber\\[1ex]
&\lesssim&\rho^{\frac{2}{s}-\frac{4}{3}}\Vert\nabla u_2\Vert_{L^{3/2}(D_{3\rho/4}(x))}+\Vert\nabla u_3\Vert_{L^{s}(D_{5\rho/8}(x))}
\nonumber\\[1ex]
&\lesssim&\rho^{\frac{2}{s}-\frac{4}{3}}\Vert\nabla u\Vert_{L^{3/2}(D_{3\rho/4}(x))}+\Vert\nabla u_3\Vert_{L^{s}(D_{3\rho/4}(x))}
\nonumber\\[1ex]
&\lesssim&\rho^{\frac{2}{s}-1}\Vert\nabla u\Vert_{L^{2,\infty}(D_{3\rho/4}(x))}+\big\Vert|z|^{-1}*\chi_{D_{\rho/4}(y)}\Delta u\big\Vert_{L^{s}(D_{3\rho/4}(x))}
\nonumber\\[1ex]
&\leq&C(\eps_0)\Big[\rho^{\frac{2}{s}-1}\Vert\nabla u\Vert_{L^{2,\infty}(D_{\rho}(x))}+
\rho^{\frac{2}{s}-\frac{2}{p}+1}\Vert F\Vert_{L^p(D_\rho(x))}\Big]\:.
\end{eqnarray}
As seen above, we can choose any $s<2/(2-p)$. 

\medskip

Suppose next that $p\in[2,\infty)$. Let $s\in(2,\infty)$ be arbitrary. Choose $0<\eps<2/s$. Then, setting $q=2-\eps$, we have
$$
\Vert F\Vert_{L^{q}(D_\rho(x))}\;\lesssim\;\rho^{\frac{2}{q}-\frac{2}{p}}\Vert F\Vert_{L^{p}(D_\rho(x))}\:.
$$
Since $s<2/(2-q)$, we have per the above discussion that
\begin{eqnarray*}
\Vert\nabla u\Vert_{L^s(D_{5\rho/8}(x))}&\leq&C(\eps_0)\Big[\rho^{\frac{2}{s}-1}\Vert\nabla u\Vert_{L^{2,\infty}(D_{\rho}(x))}+\rho^{\frac{2}{s}-\frac{2}{q}+1}\Vert F\Vert_{L^q(D_\rho(x))}\Big]\nonumber\\[1ex]
&\leq&C(\eps_0)\Big[\rho^{\frac{2}{s}-1}\Vert\nabla u\Vert_{L^{2,\infty}(D_{\rho}(x))}+\rho^{\frac{2}{s}-\frac{2}{p}+1}\Vert F\Vert_{L^p(D_\rho(x))}\Big]\:.
\end{eqnarray*}
In other words, (\ref{bell1}) holds for all $p\in(1,\infty)$, with any $s<2/(2-p)$ if $p\in(1,2)$ and any $s<\infty$ if $p\ge2$. We combine these facts by writing the (\ref{bell1}) holds for all $s<2/(2-p)_+$.

\eject


\bigskip


\begin{thebibliography}{99}
\bibitem[Ada]{Ad} Adams, Robert ``A note on Riesz potentials." Duke Math. J. 42 (1975), no. 4, 765--778.
\bibitem[Ber1]{Ber1} Bernard, Y. ``Analysis of constrained Willmore surfaces." Comm. PDE {41} (2016), no. 10, 1513--1552.
\bibitem[Ber2]{Ber2} Bernard, Y. ``Noether's theorem and the Willmore functional." Adv. Calc. Var. {9} (2016), no. 3, 217--234.
\bibitem[Ber2b]{Ber2b} Bernard, Y. ``Noether's theorem and the Willmore functional." arXiv preprint 1409.6894.
\bibitem[BR1]{BR0} Bernard, Y.; Rivi\`ere, T. ``Local Palais-Smale sequences for the Willmore functional." Comm. Anal. Geom. 19 (2011), 563--599.
\bibitem[BR2]{BR1} Bernard, Y.; Rivi\`ere, T. ``Singularity removability at branch points for Willmore surfaces." Pacific J. Math. 265 (2013), 257--311.
\bibitem[BR3]{BR2} Bernard, Y.; Rivi\`ere, T. ``Energy quantization for Willmore surfaces and applications." Ann. of Math. 180 (2014), 87--136. 
\bibitem[BR4]{BR3} Bernard, Y.; Rivi\`ere, T. ``Ends of immersed minimal and Willmore surfaces in asymptotically flat spaces." arXiv preprint 1508.01391.
\bibitem[BWW]{BWW} Bernard, Y.; Wheeler, G.; Wheeler, V.-M.``Rigidity and stability of spheres in the Helfrich model." Interfaces Free Bound. 19 (2017), 495--523.
\bibitem[BWW2]{BWW2} Bernard, Y.; Wheeler, G.; Wheeler, V.-M.``Concentration-compactness and finite-time singularities for Chen's flow" arXiv preprint 1706.01707.
\bibitem[BPP]{BPP} Bohle, C.; Peters, G.P.; Pinkall, U. ``Constrained Willmore surfaces." Calc. Var. PDE 32 (2008), 263--277.
\bibitem[Can]{Can} Canham, P.B. ``The minimum energy of bending as a possible explanation of the biconcave shape of the human red blood cell." J. Theor. Biol. 26 (1970), 61--81.
\bibitem[Che1]{BYC1} Chen, B.-Y. ``Some open problems and conjectures on submanifolds of finite type." Soochow J. Math. 17 (1991), 169--188.
\bibitem[Che2]{BYC2} Chen, B.-Y. ``Recent developments of biharmonic conjecture and modified biharmonic conjectures." arXiv preprint 1307.0245.
\bibitem[Dim1]{Dim1} Dimitri\'c, I. ``Submanifolds of $\mathbb{E}^n$ with harmonic mean curvature vector." Bull. Inst. Math. Acad. Sinica 20 (1992), 53--65.
\bibitem[Dim2]{Dim2} Dimitri\'c, I. ``Quadric representation and submanifolds of finite type." Ph.D. Thesis, Michigan State University (1989).
\bibitem[Gia]{Gia} Giaquinta, M. ``Multiple Integrals in the Calculus of Variations and Nonlinear Elliptic Systems." PUP, Princeton (1983).
\bibitem[HV]{HV} Hasanis, T.; Vlachos, T. ``Hypersurfaces in $\mathbb{E}^4$ with harmonic mean curvature vector field." Math. Nachr. 172 (1995), 145--169.
\bibitem[Hel]{Hel} H\'elein, F. ``Harmonic Maps, Conservation Laws, and Moving Frames." Cambridge Tracts in Mathematics, 150. Cambridge University Press (2002).
\bibitem[Hef]{Hef} Helfrich, W. ``Elastic properties of lipid bilayers: theory and possible experiments."  Z. Naturforsch., C28 (1973), 693--703. 
\bibitem[Hub]{Hub} Huber, Alfred ``On subharmonic functions and differential geometry in the large." Comment. Math. Helv. 32 (1957), 181--206.
\bibitem[Jia]{Jia} Jiang, G.Y. ``2-harmonic isometric immersions between Riemannian manifolds." Chinese Ann. Math. Ser. A, 7 (1986), 130--144.
\bibitem[KMR]{KMR} Keller, L.G.A.; Mondino, A.; Rivi\`ere, T. ``Embedded surfaces of arbitrary genus minimizing the Willmore energy under isoperimetric constraint." Arch. Rat. Mech. Anal. 212 (2014), 645--682.
\bibitem[KS1]{KS1} Kuwert, E.; Sch\"atzle, R. ``The Willmore flow with small initial energy." J. Diff. Geom. 57 (2001), 409--441.
\bibitem[KS2]{KS2} Kuwert, E.; Sch\"atzle R. ``Removability of point singularities of Willmore surfaces." Ann. of Math. 160 (2004), no. 1, 315--357.
\bibitem[KS3]{KS3} Kuwert, E.; Sch\"atzle, R. ``Minimizers of the Willmore energy under fixed conformal class." J. Diff. Geom. 93 (2013), 471--530.
\bibitem[Lnk]{Lnk} Link, F. ``Gradient flow for the Willmore functional in Riemannian manifolds of bounded geometry.'' arXiv preprint 1308.6055.
\bibitem[MN]{MN} Marques, F.C.; Neves, A. ``Min-max theory and the Willmore conjecture.'' Ann. of Math. 179 (2014), 683--782.
\bibitem[MPW]{MPW} McCoy, J; Parkins, S.; Wheeler, G. ``The geometric triharmonic heat flow of immersed surfaces near spheres.'' Nonlinear Anal. 161 (2017), 44--86.
\bibitem[MW1]{MW1} McCoy, J; Wheeler, G. ``A classification theorem for Helfrich surfaces.'' Math. Ann. 357:4 (2013), 1485--1508.
\bibitem[MW2]{MW2} McCoy, J; Wheeler, G. ``Finite-time singularities for the locally constrained Willmore flow of surfaces.'' Comm. Anal. Geom. 24 (2016), no. 4, 843--886.
\bibitem[MWW]{MWW} Metzger, J; Wheeler, G.; Wheeler, V.-M.. ``Willmore flow of surfaces in Riemannian spaces I: Concentration-compactness.'' arXiv preprint 1308.6024.
\bibitem[MS]{MS}  M\"uller, Stefan ;  \v{S}ver\'ak, Vladim\'ir ``On surfaces of finite total curvature." J. Diff. Geom. 42 (1995), no. 2, 229--258.
\bibitem[Riv1]{Riv1} Rivi\`ere, T. ``Analysis aspects of the Willmore functional." Invent. Math. 174 (2008), no. 1, 1--45.
\bibitem[Riv2]{Riv2} Rivi\`ere, T.  ``Variational principles for immersed surfaces with $L^2$-bounded second fundamental form.'' J. reine angew. Math. (2013).
\bibitem[Riv3]{ParksCity} Rivi\`ere, T. ``Weak immersions of surfaces with $L^2$-bounded second fundamental form''. Parks City lecture notes, July 2013. 
\bibitem[Sch]{Sch} Sch\"atzle, R. ``Conformally constrained Willmore immersions." Adv. Calc. Var. 6 (2013), 375--390.
\bibitem[Scy]{Scy} Schygulla, J. ``Willmore minimizers with prescribed isoperimetric ratio'', Arch. Rat. Mech. Anal. 203 (2012), no. 3, 901--941.
\bibitem[Tar]{Tar} Tartar, Luc ``An Introduction to Sobolev Spaces and Interpolation Spaces." Lectures notes of the Unione Matematica Italiana, no. 3 (2007).
\bibitem[Wei]{Wei} Weiner, J. ``On a problem of Chen, Willmore, et al.'' Indiana Univ. Math. J. 27 (1978), no. 1, 19--35.
\bibitem[Whe]{Whe} Wheeler, G. ``Fourth order geometric evolution equations'' Bull. Aust. Math. Soc. 82 (2010). PhD thesis.
\bibitem[Whe1]{Whe1} Wheeler, G. ``Gap phenomena for a class of fourth-order geometric differential operators on surfaces with boundary.'' Proc. AMS 143 (2015), no. 4, 1719--1737.
\bibitem[Whe2]{Whe2} Wheeler, G. ``Surface diffusion flow near spheres.'' Calc. Var. PDE 44 (2012), no 1, 131--151.
\bibitem[Whe3]{Whe3} Wheeler, G. ``Chen's conjecture and $\varepsilon$-superbiharmonic submanifolds of Riemannian manifolds.'' Internat. J. Math. 24 (2013), no. 4. 1350028.
\bibitem[Wil]{Wil} Willmore, T. J. ``Note on embedded surfaces." Ann. Stiint. Univ. ``Al. I. Cuza" Iasi. Sect. I a Mat. (N.S.) 11B (1965), 493--496.
\end{thebibliography}
\end{document}